\documentclass[a4paper,11pt]{article}
\usepackage{epsf}
\usepackage{enumitem}
\usepackage[french]{babel}
\usepackage{graphicx}
\usepackage{mathptmx}
\usepackage{makeidx}
\usepackage{amsmath}
\usepackage{amssymb}
\usepackage{amsfonts}
\usepackage[bottom]{footmisc}
\usepackage{multicol}
\usepackage{type1cm}
\usepackage{mathrsfs}
\usepackage{url}
\usepackage{color}
\usepackage[pagewise]{lineno}
\usepackage{setspace}
\usepackage{float}

\usepackage[T1]{fontenc}
\usepackage[applemac]{inputenc}
\newtheorem{theoreme}{Theorem}[section]
\newtheorem{lemme}[theoreme]{Lemme}
\newtheorem{proposition}[theoreme]{Proposition}

\newtheorem{definition}[theoreme]{Definition\rm}
\newtheorem{hypothese}[theoreme]{Hypothese\rm}
\newtheorem{remarque}{\bf Remarque}

\title{Entrée-sortie dans le halo d'une courbe lente semi-stable}
\author{C. Lobry}
\date{\today}

\newcommand{\bitbul}{\begin{itemize}[label = \textbullet]}
\newcommand{\bittiret}{\begin{itemize}[label = -]}
\newcommand{\bito}{\begin{itemize}[label =$\circ$]}
\newcommand{\bit}{\begin{itemize}}
\newcommand{\fit}{\end{itemize}}
\newcommand{\ben}{\begin{enumerate}}
\newcommand{\fen}{\end{enumerate}}

\newcommand{\fin}{\end{document}}
\newcommand{\beq}{\begin{equation}}
\newcommand{\feq}{\end{equation}}
\newcommand{\dcom}{\begin{quote}\begin{small}}
\newcommand{\fcom}{\end{small}\end{quote}}

\newcommand{\bc}{\begin{center}}
\newcommand{\fc}{\end{center}}
\newcommand{\emat}{\mathrm{e}}
\newcommand{\ch}{\mathrm{cosh}}
\newcommand{\sh}{\mathrm{sinh}}
\newcommand{\eps}{\varepsilon}
\newcommand{\Rmat}{\mathbb{R}}
\newcommand{\Nmat}{\mathbb{N}}

\newcommand{\arc}{\overset{\displaystyle  \frown}}
\def\1{{\rm 1\mskip-4.4mu l}}
\newcommand{\eset}[1]{{
  \mathchoice
    {\left\{\!\!\left\{ #1 \right\}\!\!\right\}} 
    {\left\{\!\left\{ #1 \right\}\!\right\}} 
    {} 
    {} 
  }
}

\begin{document}
\maketitle

\section*{Introduction}
Je m'intéresse au portrait de phase du système :
 \beq \label{Smintro}
 S_m\quad  \quad \left\{
\begin{array}{lcl}
\displaystyle \frac{dx}{dt}& =& 1\\[6pt]
\displaystyle  \frac{dy}{dt} &=& \displaystyle \frac{1}{\eps}\sqrt{m^2+y^2} \Big( f(x) -y\Big)
 \end{array} 
 \right.
\feq
lorsque les paramètres $m$ et $\eps$ sont petits. Ce système apparait de façon naturelle dans l'étude du phénomène d'{\em inflation}, important en dynamique des populations (voir \cite{BLSS21},\cite{HOLTPNAS20} et \cite{KAT21}). 

Il y a un peu plus de quarante ans ont été définies des solutions particulières de systèmes de la forme :
 \beq \label{eqLR2}
 \mathrm{LR}\quad \quad  \quad \left\{
\begin{array}{lcl}
\displaystyle \frac{dx}{dt}& =& F(m,x,y)\\[6pt]
\displaystyle  \frac{dy}{dt} &=& \displaystyle \frac{1}{\eps}G(m,x,y)
 \end{array} 
 \right.
\feq
\begin{minipage}{0.5 \textwidth}
 dans le cas particulier où l'ensemble  des zéros de $G(0,x,y) = 0$ est constituée de deux courbes $\gamma_i$, appelées ''courbes lentes'', qui se coupent transversalement et où $G$  change de signe à la traversée de chaque $\gamma_i$.  Dans ce cas chaque courbe $\gamma_i$ se décompose en $\gamma_i^{\;a}$ (partie attractive) et $\gamma_i^{\;r}$ (partie répulsive). Lorsque $\eps$ est petit on montre que pour certains $m$ de l'ordre de $\emat^{-\frac{1}{\eps}}$ il existe, ce qui est surprenant, des 
 \end{minipage} $\quad$
\begin{minipage}{0.45 \textwidth}
\includegraphics[width=1\textwidth]{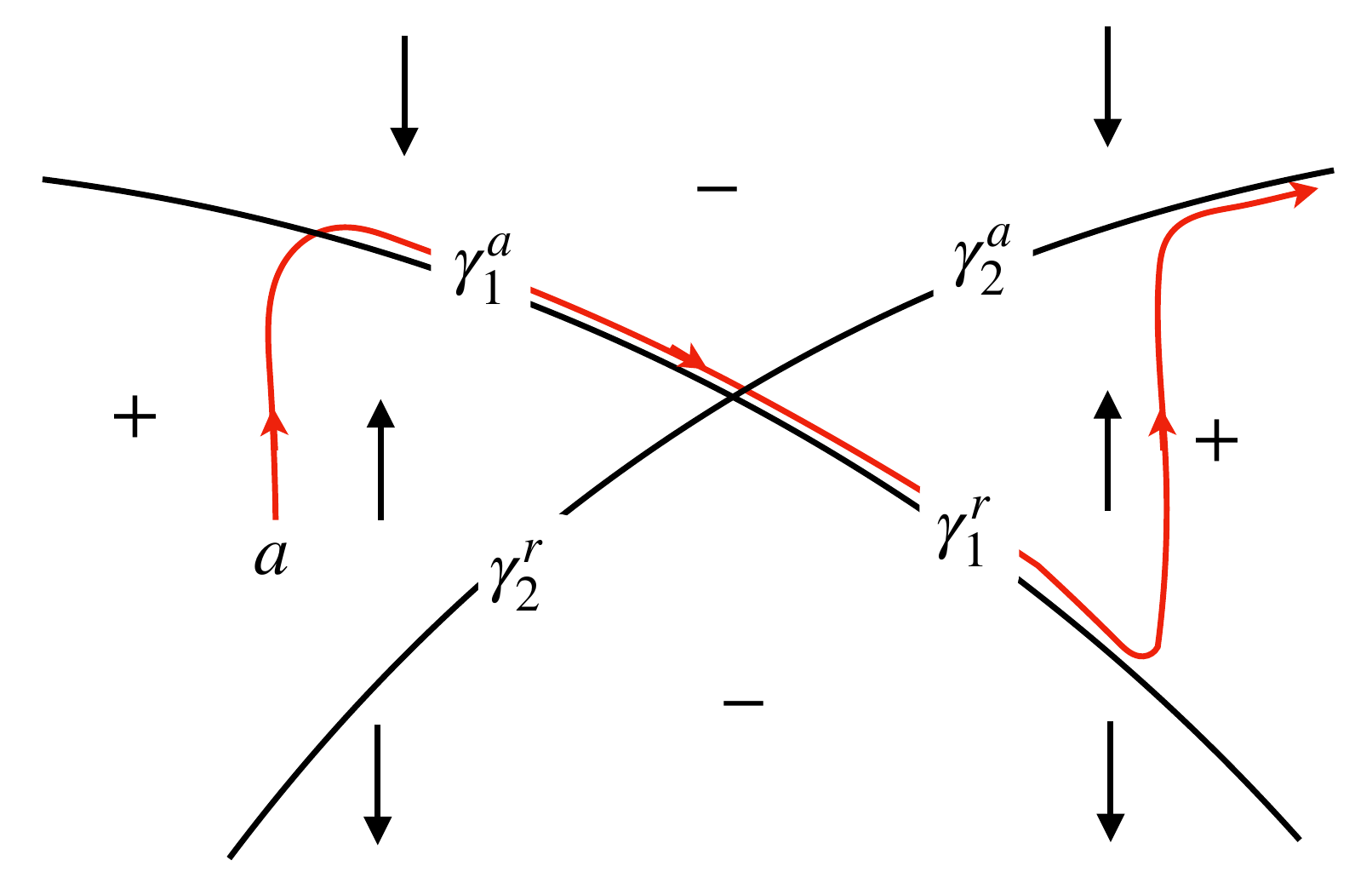}
\end{minipage} \\[2pt]
 solutions qui, comme  la solution de la figure ci-contre,  après avoir longé la partie attractive $\gamma_1^{\;a}$ de la courbe lente $\gamma_1$ continuent à longer, pendant une durée significative, la partie répulsive $\gamma_1^{\;r}$. Ces solutions ont été appelées des {\em canards} par leurs inventeurs : \cite{BEN81, BCDD81}. On trouvera en ligne sur Scholarpedia l'article \cite{WEC07} qui fait l'historique de cette découverte et de ses développements. Le lecteur intéressé trouvera également dans \cite{FRU09} une revue des applications de la théorie des {\em canards} à la théorie des bifurcations dynamiques.

Le système qui m'intéresse ressemble à celui dont nous venons de parler mais pas exactement. Pour $m=0$ on a :
 \beq 
 S_0\quad  \quad \left\{
\begin{array}{lcl}
\displaystyle \frac{dx}{dt}& =& 1\\[6pt]
\displaystyle  \frac{dy}{dt} &=& \displaystyle \frac{1}{\eps}|y| \Big( f(x) -y\Big)
 \end{array} 
 \right.
\feq
Ici, c'est $ \displaystyle \frac{1}{\eps}|y| \Big( f(x) -y\Big)$ qui correspond à la fonction $G$, les courbes lentes $\gamma_1$ et $\gamma_2$ sont respectivement le graphe de $f$ et l'axe  $y = 0$ mais la distribution des  signes \\[2pt]
\begin{minipage}{0.45 \textwidth}
 est différente.
 Ainsi la courbe lente $\gamma_1$ ne passe plus comme précédemment d'attractive à répulsive mais reste attractive alors que $\gamma_2$ (ici l'axe $y = 0$ ) est semi attractif, semi répulsif. Nous allons montrer l'existence, toujours pour des valeurs exponentiellement petites de $m$  par rapport à $\eps$, de solutions qui, comme celle de la figure, longent pendant des durées significatives des portions de $\gamma_2$, donc des portions  {\em non stables} de la courbe lente.
\end{minipage}$\quad $
\begin{minipage}{0.48 \textwidth}
\includegraphics[width=1\textwidth]{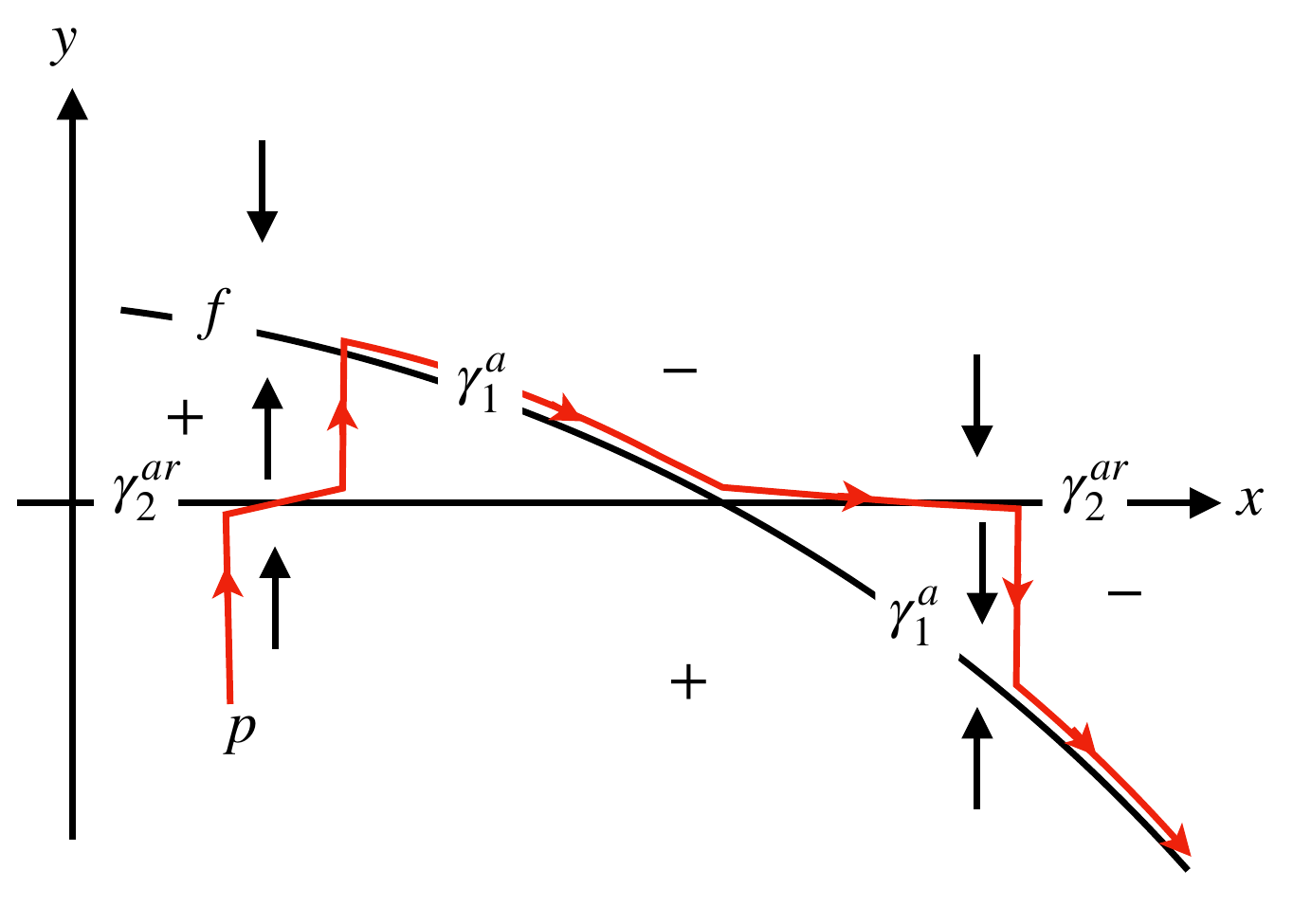}
\end{minipage} \\[6pt]
 Il n'est plus question de ''canards'' au sens de la définition usuelle mais le phénomène est de même nature comme je vais le montrer.
 
 Pour démontrer l'existence de ces solutions j'utiliserai les méthodes des inventeurs des ''canards'', c'est à dire l'Analyse Non Standard (ANS) dans le formalisme I.S.T. de Nelson et les techniques de changements de variables développées notamment dans \cite{BEN81, BCDD81}. 
 Le lecteur familier des travaux des inventeurs des  ''canards'' \cite{BEN81, BCDD81} constatera que cet article   est un simple exercice application des méthodes qui y sont introduites. Pour que ce texte soit compréhensible par un lecteur non familier des techniques de l'ANS j'ai essayé de réduire au minimum possible les concepts et la terminologie liée à l'ANS et les ai décrits dans une courte annexe.  
 
  L'objectif est de dresser le portrait de phase de \eqref{Smintro} et, en raison des applications possibles, déjà évoquées, en dynamique des populations, je traite  le cas ou $f$ n'est pas nécessairement continue.
  Dans la première section j'énonce et illustre sur des exemples  le résultat principal  de cet article : la théorème d'approximation des solutions. 
  Dans la section suivante je démontre le théorème puis  je termine par une application : la démonstration d'une conjecture énoncée par G. Katriel dans \cite{KAT21}. Les démonstrations sont géométriques et se lisent essentiellement sur des figures qui, si l'article est lu sur une double page, se trouvent dans toute la mesure du possible en regard du texte. 
\newpage 
 
\tableofcontents
\newpage
\section{Le système  $S_m$}
\subsection{Vocabulaire et notations} 
Il est nécessaire d'introduire tout de suite un vocabulaire ANS  minimal et les notations associées. Le lecteur qui ne connait pas les méthode de l'ANS peut se contenter du sens intuitif des mots introduits mais il trouvera, s'il le désire, des informations plus complètes en annexe.

En analyse non standard tout tourne autour du fait que la phrase : 
 \bitbul
 \item $\eps$ est un nombre réel strictement positif infiniment petit (i.p.) 
 \fit
 a un sens formel précis (voir annexe \ref{ANS}). A partir de là on définit immédiatement :
 \bitbul
 \item $x \in \Rmat$ est infiniment grand (i.g.) $\stackrel{df}{\longleftrightarrow}  \exists \, \eps>0\; \mathrm{i.p.\;tel\;que}\; x > \frac{1}{\eps}$ 
 \item $x \in \Rmat$ est limité si $\stackrel{df}{\longleftrightarrow}  |x|$ n'est pas infiniment grand.
 \fit
 et on introduit les notations :
 \bitbul
 \item $x \sim y \stackrel{df}{\longleftrightarrow}  |x-y|$ i.p.
 
 \item $x \lnsim y \stackrel{df}{\longleftrightarrow} x<y|$ et $ |y-x] \nsim 0$
 
\item Le {\em halo} d'un sous ensemble $E$  de $\Rmat^2$ est constitué des points qui sont infiniment proches des points de $E$.

\item Les réels infiniment grands et infiniment petits sont {\em nonstandard} mais il y en a beaucoup d'autres comme, par exemple, $1+\eps$ lorsque $\eps$ est infiniment petit\footnote{
Les nombres en écriture décimale avec un nombre fini de décimales fournissent une  image de la droite réelle nonstandard : on peut décider que des nombres avec 6 décimales sont standard et que des nombres avec 12 décimales sont non standard. Un infiniment petit est un nombre de la forme $0,000\;000\;***\;***$.}. Tout réel limité est infiniment proche d'un unique réel standard.

\item Une fonction $f$ est dite {\em limitée} si il existe $m$ limité tel que $\forall x \;|f(x)|< M$ et  $C^1$-{\em limitée } si elle est limitée ainsi que sa dérivée.
 \fit
\subsection{Les équations}\label{equationsdebase}
Dans tout cet article on s'intéresse au système :
 \beq \label{Sm}
 S_m\quad  \quad \left\{
\begin{array}{lcl}
\displaystyle \frac{dx}{dt}& =& 1\\[6pt]
\displaystyle  \frac{dy}{dt} &=& \displaystyle \frac{1}{\eps}\sqrt{m^2+y^2} \Big( f(x) -y\Big)
 \end{array} 
 \right.
\feq
\paragraph{Hypothèses et notations sur $f$ :}\label{hyp}
On suppose que $f$ est $C^1$-limitée par morceaux, c'est à dire :
\ben
\item Il existe une suite discrète $D = \{x_n \;;\; n\in \mathbb{Z^*}\}$ 

\item Sur chaque intervalle $[x_n,x_{n+1}[$ la fonction $f$ est la restriction d'une fonction  $f_n$, $C^1$-limitée, définie sur $\Rmat$ tout entier telle que $f_n(x) = 0 \Rightarrow f_n'(x) ≠0$ (donc les zéros de $f_n$ sont isolés).

\item On note $\mathcal{C}$ le graphe de $f$ et $\theta_n \,:\, n \in J$ ($\theta_n < \theta_{n+1}$) la suite discrete des valeurs de changement de signe de $f$. Pour un $x$ donné on note $\theta(x)$ le premier $\theta_n$ plus grand que $x$.
\fen

 Sous ces hypothèses\footnote{
Ces hypothèses ne sont  pas absolument nécessaires  mais permettent d'alléger les écritures notamment  en entrainant que les solutions ne peuvent pas tendre vers l'infini pour des valeurs finies de $t$.}, pour toute condition initiale, le système $S$ admet une solution unique définie pour tout $t$. Précisément si $x_0 \in [x_n,x_{n+1}[$ on intègre $S_m$, avec $f_n$ à la place de $f$, à partir de la condition initiale $(x_0,y_0)$ ; comme $f_n$ est bornée, la solution est définie jusqu'à ce qu'elle rencontre la verticale $x = x_{n+1}$ en un point $(x_{n+1}, y_{n+1})$ qui sert de nouvelle condition initiale et ainsi de suite.
\bitbul
\item On note $(x(t, (x_0,y_0),m),y(t,(x_0,y_0),m))$ la solution de  $S_m$ de condition initiale $(x_0,y_0)$ à l'instant $0$. On a  $x(t, (x_0,y_0),m) = x_0+t$.
\fit
Il faut bien noter que ce théorème d'existence de solutions est vrai pour tout $\eps >0$ et tout $m$, qu'ils soient standard ou non. A partir de maintenant :

\bitbul
\item On suppose que $\eps >0$ { infiniment petit} est donné une fois pour toutes et 
$$m= \emat^{\frac{\rho}{\eps}}\quad \rho < 0$$
 est un paramètre. 
\fit
On s'intéresse à l'évolution  portrait de phase  de $S_m$ en fonction de $\rho$.

Notons que ce n'est pas le souci de généralité maximum qui conduit à considérer des fonctions $f$ continues par morceaux plutôt que simplement continues mais des contextes où des modélisations avec des fonctions discontinues sont plus naturelles (voir \cite{BLSS21}). De plus cette plus grande généralité n'introduit pas de difficulté supplémentaire dans le type de preuves que j'utilise.
$ \quad$

\subsection{Le système contraint et le théorème d'approximation.}\label{contraint}
\begin{figure}[t]
  \begin{center}
 \includegraphics[width=0.9\textwidth]{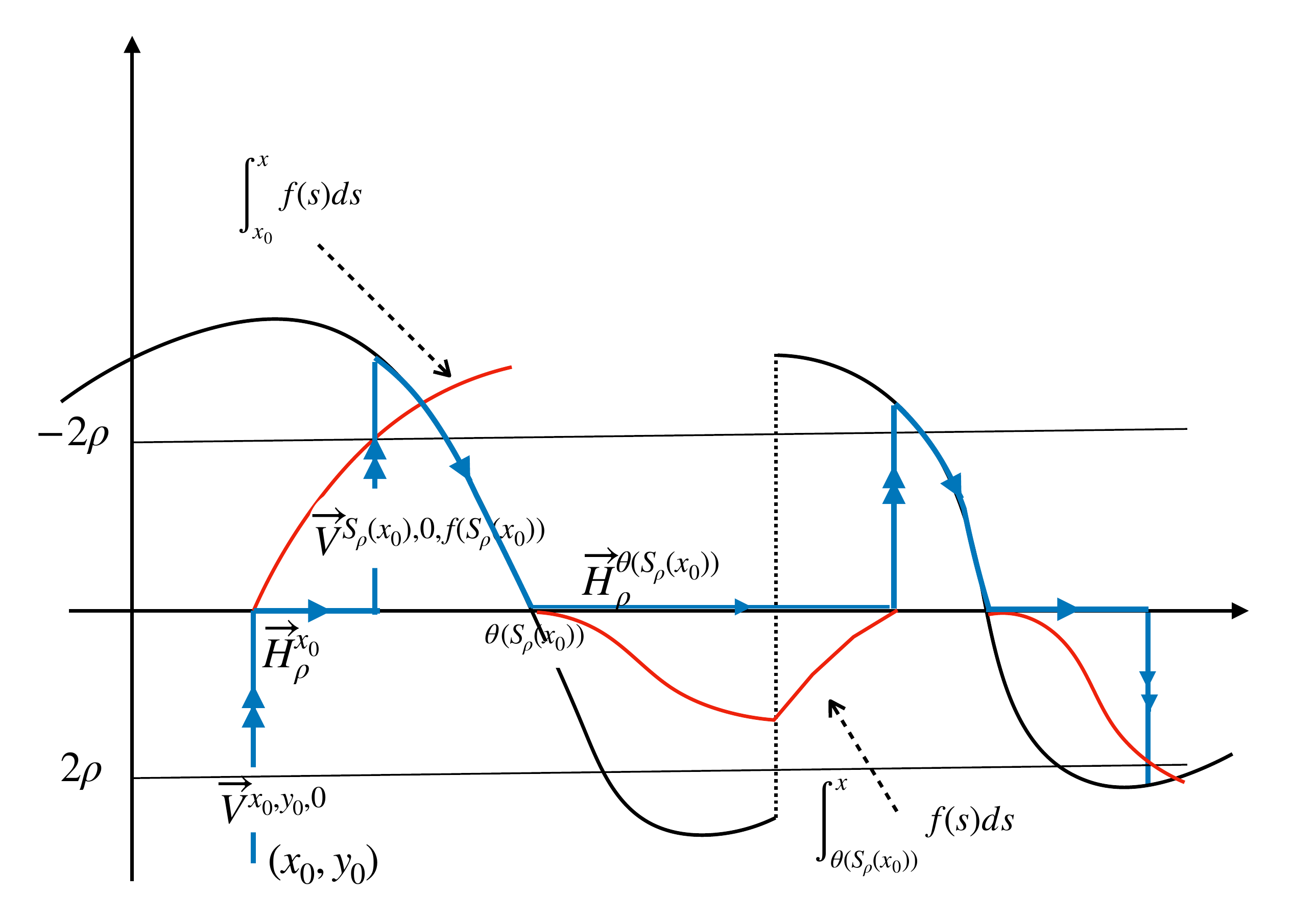}
 \caption{C-trajectoire de \eqref{SCm} issue de $(x_0,y_0)$. ; explications sous-section \ref{contraint} } \label{figSCm}
 \end{center}
 \end{figure}

Au système $S_m$ on associe le "système contraint" :
 \beq \label{SCm}
 S_m^0\quad  \quad \left\{
\begin{array}{lcl}
\displaystyle \frac{dx}{dt}& =& 1\\[6pt]
\displaystyle 0 &=& \displaystyle \sqrt{m^2+y^2} \Big( f(x) -y\Big)
 \end{array} 
 \right.
\feq
Ce n'est pas un système différentiel mais, toutefois, on peut lui associer des ''pseudo-trajectoires'' de la façon suivante\footnote{
La définition de notion de solution pour des "systèmes contraints", de la forme ci-dessus, mais en dimension plus grande, a fait l'objet d'importantes études (par exemple \cite{TAK76}) mais il n'est pas nécessaire d'y faire référence ici car elles ne portent pas sur le lien entre les solutions de $S_m$ et de  $S^0_m$. En revanche on pourra consulter \cite{WEC07} et ses références  sur l'existence de {\em canards} en dimension plus grande que 2.}.
\begin{definition} On se donne $\rho <0$. 0n appelle :\\
\bitbul
\item \textbf{Segment vertical}   un segment de la forme $$ \overrightarrow{V}^{x,[a,b]}=[(x,a), (x,b)]$$
tel que $0 \not \in ]a,b[$ et $b = 0$ ou $b = f(a)$ orienté vers le haut en dessous du graphe de $f$ et vers le bas en dessus. 
	
\item \textbf{Segment de courbe lente } une partie du graphe   de $f$ de la forme :
$$ \mathcal{C}^x = \{(s,f(s))\,:\, x \leq s \leq \theta (x)\}$$ 
où $\theta(x)$ est le premier changement de signe de $f$ qui suit $x$ (cf. hypothèses \ref{hyp}.3.), orienté vers la droite.

\item \textbf{Segment horizontal}  un segment orienté vers la droite de la forme $\overrightarrow{H}_{\rho}^{x} = [(x, 0),(0,S_{\rho}(x))]$ où $S_{\rho}(x)$ est défini comme le plus petit   $x^*>x$ (éventuellement $+\infty$) tel que :
\beq
\displaystyle \int_x^{x^*} f(s)ds \in \{+2\rho,\,0, -2\rho\}\feq
L'extrémité  $(S_{\rho}(x),0)$ du segment horizontal  s'appelle ''point de sortie''.
\fit
\end{definition}
\begin{proposition} \label{rderho} 
Pour $x$ fixé, $S_{\rho}(x)$ est une fonction décroissante de $\rho $ qui tend vers $x$ lorsque $\rho$ tend vers $0$.
\end{proposition}
\textbf{Preuve.} Evident. $\Box$
\begin{definition}
Une { \em C-trajectoire} (pour ''trajectoire du système contraint'') de (\ref{SCm}) est une suite $\Gamma_i$ de segments mis bout à bout (i.e. l'origine de $\Gamma_{i+1}$ est l'extrémité de $\Gamma_i$) telle que :
\bito

\item à un segment vertical $ \overrightarrow{V}^{x,[a,f(x)]}$ succède le segment de courbe lente $ \mathcal{C}^x$ 
\item à un segment vertical $ \overrightarrow{V}^{x,[a,0]}$ succède le segment horizontal $\overrightarrow{H}_{\rho}^{x}$ 
\item à un segment horizontal $\overrightarrow{H}_{\rho}^{x}$ succède le segment vertical $ \overrightarrow{V}^{S_{\rho}(x),[0,f(_{\rho}(x))]}$
\item  à un segment de courbe lente  $ \mathcal{C}^x$ succède le segment  horizontal $\overrightarrow{H}_{\rho}^{\theta(x)}$
\fit
Ces règles permettent d'associer à toute ''condition initiale'' $(x_0,y_0)$ une unique C-trajectoire comme on le constate sur la figure \ref{figSCm}.
\end{definition} 
\begin{definition} Soit $\Gamma$ un arc   $s \mapsto (\alpha (s),\beta(s))$, ($s \in [a,b]$), de $\Rmat^2$ et $t \mapsto (x(t),y(t))$, ($ t\in [t_1,t_2]$), un autre arc de $\Rmat^2$. On dit que 
$ (x(t),y(t)) \;\mathrm{longe}\; \Gamma$
si il existe un paramétrage $t \mapsto s(t)$, ($t \in [t_1,t_2]$),  de $\Gamma$ tel que, 
pour tout $t$, $d((x(t),y(t)),(\alpha(s(t),\beta(s(t))) \sim 0$, où $d$ est la distance naturelle de $\Rmat^2$.
\end{definition}
On peut maintenant énoncer le résultat principal de cet article.

\begin{theoreme} \label{theoreme} \textbf{Théorème d'approximation.}
Soit $(x_0,y_0)$ une condition initiale limitée, telle que $y_0 \not \sim 0$, et $t \mapsto (x(t),y(t))$ la solution du système différentiel \eqref{Sm} qui en est issue à l'instant $t_0$. Alors $(x(t),y(t))$ est infiniment proche de la C-trajectoire  de \eqref{SCm} issue du même point. Plus précisément, si $\Gamma_1, \Gamma_2,\cdots, \Gamma_n,\cdots$ sont les segments successifs de la C-trajectoire, il existe une suite d'instants, $t_1,t_2,\cdots,t_n,\cdots$ tels que sur $[t_{n-1},t_n]$ la trajectoire $(x(t),y(t))$ longe le segment $\Gamma_n$.
\end{theoreme}

La dénomination ''point de sortie'' pour l'extrémité $(S_{\rho}(x),0)$ du segment horizontal $\overrightarrow{H}_{\rho}^{x} $ vient de ce que la trajectoire qui avait ''pénétré'' au point $(x,0)$ dans le halo de $y = 0$ en ''ressort'' au point $(S_{\rho}(x),0)$. Pour cette raison on appelle ''retard'' (sous entendu ''à la sortie'') la quantité $R_{\rho}(x) = S_{\rho}(x) - x$.\\

\noindent\textbf{ L'hypothèse $y_0 \not \sim 0$ est essentielle comme nous le verrons au paragraphe \ref{ci0}}

\paragraph{Remarque} {\em Pour ceux qui connaissent I.S.T.).} En utilisant toute la puissance de I.S.T. on peut introduite le concept d'{\em ombre} d'un ensemble (le {\em standardisé} du halo). Le discours précédent est reformulation  de la proposition : {\em L'ombre de la trajectoire issue du point $a$ est la C-trajectoire issue de $a$}. En se limitant à une version simplifiée de I.S.T. (comme je le fais ici) on se prive d'outils facilitant la rédaction comme le serait  l'usage de l'ombre mais on diminue le prix à payer pour la pratique de l'ANS. Je discute ce point dans l'annexe \ref{cheapANS}

\newpage
\subsection{Illustrations du théorème d'approximation.}\label{illustrations}

\begin{figure}
  \begin{center}
 \includegraphics[width=0.9\textwidth]{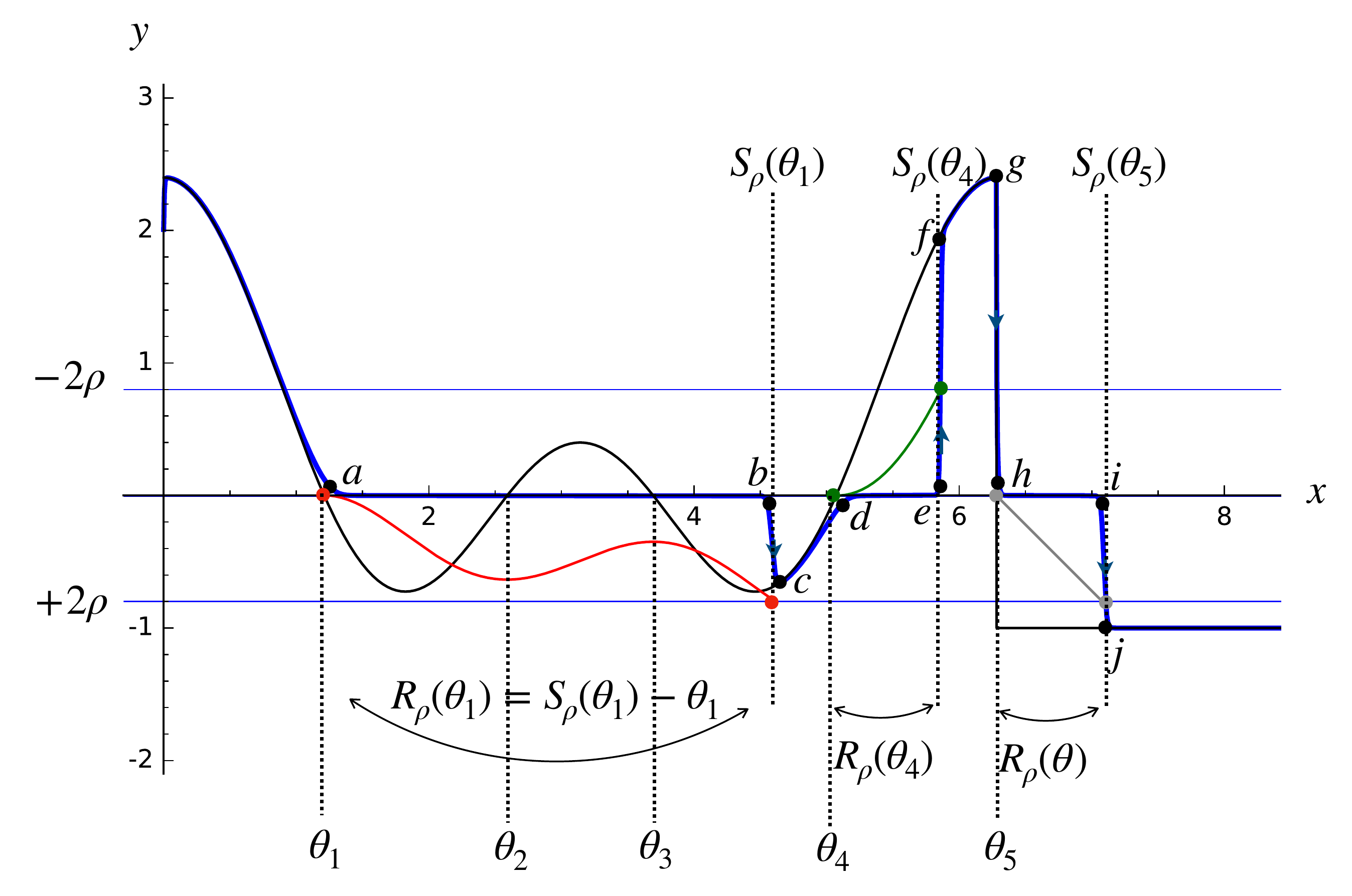}
 \caption{Prédictions du théorème d'approximation. Explications section \ref{illustrations}.} \label{retard5}
 \end{center}
 \end{figure}

 \paragraph{Commentaires sur la figure \ref{retard5}. } 
 
J'ai simulé le système $S_m$ avec :
 \bitbul
 \item La fonction $f$ définie par :
 \beq
 \left\{\;
\begin{array}{lcl}
\displaystyle x \leq 2\pi& \Rightarrow & f(x) = \cos(x)+\cos(2x)+0.4 \\
\displaystyle x> 2\pi& \Rightarrow & f(x) = -1
\end{array}
\right.
\feq
Son graphe est tracé en noir. Elle est nulle pour les valeurs $\theta_1,\theta_2,\theta_3$ et $\theta_4$ qu'il est inutile de préciser et elle est discontinue en $\theta_5 = 2\pi$.

\item $\eps = 0.01\quad \quad  \rho = -0.4\quad \quad \displaystyle m =  \emat^{-\frac{0.4}{\eps}}$
 \fit
 En bleu épais on a simulé la trajectoire $(x(t),y(t))$ issue du point $(0,2)$ à l'instant $t_0 = 0$. On voit, qu'en dehors des ''angles un peu arrondis'' elle présente une suite de ''segments'' comme la C-trajectoire. Nous détaillons cette suite :
 \ben
 \item Il existe $t_1 \sim 0$ tel que sur $[t_0,t_1]$ la solution longe le segment vertical $ \overrightarrow{V}^{0,[2, 2.4]}$.
 \item Soit $\theta_1$ la première valeur après $x = 0$ pour laquelle  $f$ change de signe ; il existe $t_2 \sim \theta_1$ tel que sur $[t_1,t_2]$ la solution longe le segment de courbe lente $\mathcal{C}^0 = \{ (x,f(x)\;,\; x \in[0,\theta_1]\}$ jusqu'au point $a \sim(\theta_1,0)$. Le point $a$ est un point d'entrée dans le halo de la droite $y = 0$.
 \item Détermination du point de sortie   suivant $(\theta_1,0)$. En rouge on a tracé le graphe de la fonction :
 $$x \mapsto \int_{\theta_1}^xf(s)ds$$
 qui rencontre la droite $y = 2\rho$ au point d'abscisse $S_{\rho}(\theta_1) = \theta_1+R_{\rho}(\theta_1)$ ; Il existe $t_3 \sim S_{\rho}(\theta_1)$ tel que sur $[t_2,t_3]$ la trajectoire longe le segment horizontal $[(\theta_1,0),(S_{\rho}(\theta_1),0)]$ jusqu'au point $b \sim (S_{\rho}(\theta_1),0)$. On remarquera qu'il n'y a aucune raison pour que $S_{\rho}(\theta_1)$ précède $\theta_2$ ou $\theta_3$.
 \item Il existe $t_4 \sim S_{\rho}(\theta_1)$ tel que sur $[t_3,t_4]$ la solution longe le segment vertical $[(S_{\rho}(\theta_1),0),(S_{\rho}(\theta_1),f(S_{\rho}(\theta_1))]$ jusqu'au point $c \sim (S_{\rho}(\theta_1), f(S_{\rho}(\theta_1))$
 \item Soit $\theta_4$ la première valeur après $x = S_{\rho}(\theta_1))$ pour laquelle  $f$ change de signe ; il existe $t_5 \sim \theta_4$ tel que sur $[t_4,t_5]$ la solution longe le segment de graphe $\mathcal{C}^{S_{\rho}(\theta_1)} = \{ (x,f(x)\;,\; x \in[S_{\rho}(\theta_1),\theta_4]$ jusqu'au point $d \sim(\theta_4,0)$. Le point $d$ est un point d'entrée dans le halo de la droite $y = 0$.
 \item Détermination du point de sortie suivant $S_{\rho}(\theta_4,0)$. En vert on a tracé le graphe de :
 $$x \mapsto \int_{\theta_4}^xf(s)ds$$
 dont l'intersection avec la droite $y = -2\rho$ détermine le point de sortie $(S_{\rho}(\theta_4),0)$. Il existe $t_6\sim S_{\rho}(\theta_4)$ tel que sur $[t_5,t_6]$ la solution longe le segment horizontal $[(\theta_4,0),(S_{\rho}(\theta_4),0)]$.
 \item Je laisse au lecteur le soin d'interpréter les segments suivants de trajectoire : $ef,\; fg,\; gh,\;hi,\; ij$
 \fen 
\noindent \textbf{Remarque}. 
 Les prédictions du théorème d'approximation sont très précises (on peut le vérifier en agrandissant la figure \ref{retard5}) alors que les quantités manipulées dans le schéma numérique sont remarquablement petites comme par exemple $m = 4. 24\, 10^{-18}$ ; ceci est dû au fait que la représentation des réels en notation "virgule flottante" permet de représenter de très petits nombres ($10^{-250}$ avec mon logiciel). Si l'on fait le changement de variable $z = y-1$ qui translate la figure autour de la droite $z = 1$ les prédictions du théorème d'approximation seront beaucoup moins efficaces ; ce point est analysé dans \cite{LOB92}.
 
\begin{figure}[t]
  \begin{center}
 \includegraphics[width=1\textwidth]{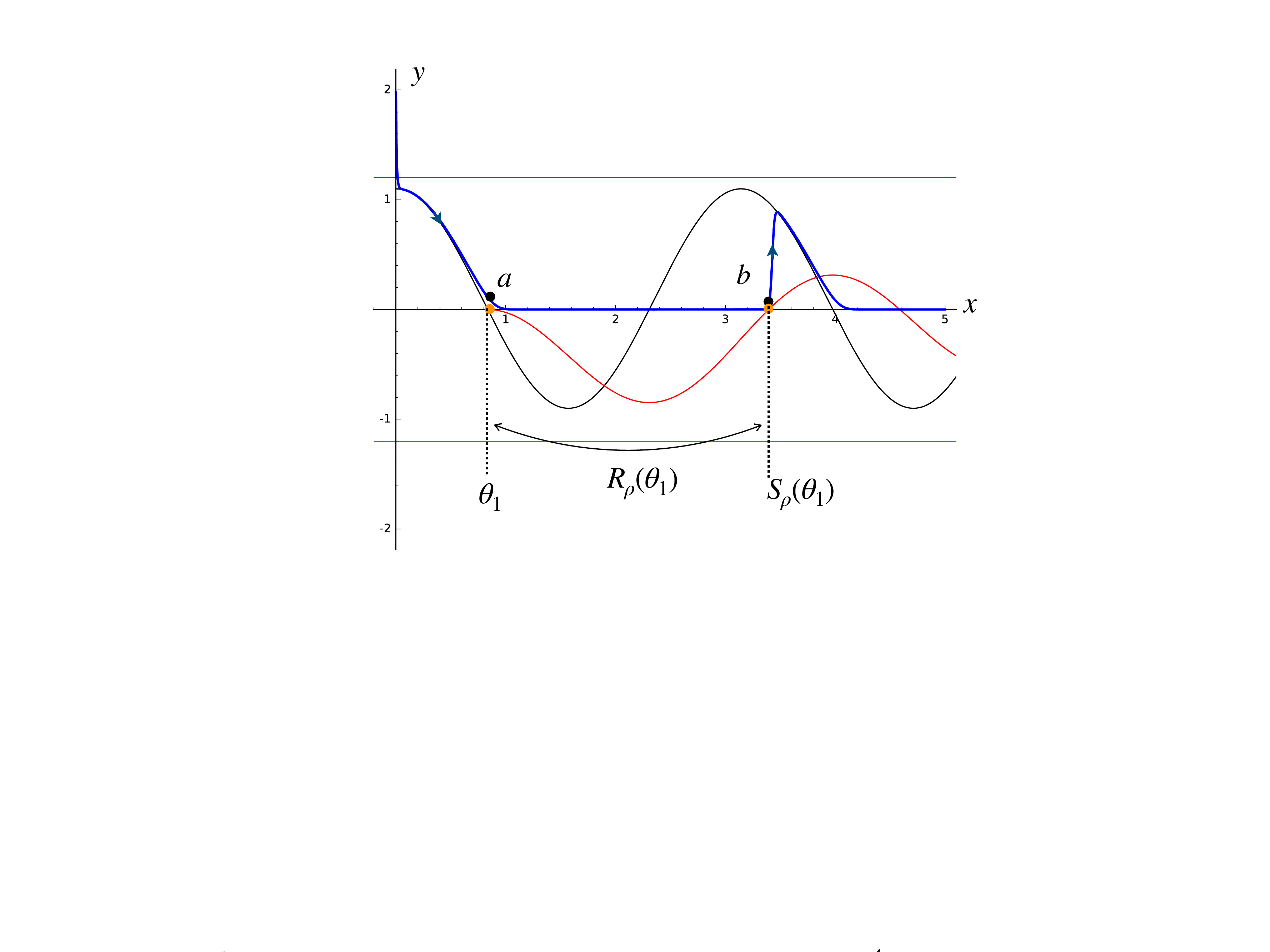}
 \caption{La trajectoire ne travese pas l'axe. Explications dans le texte } \label{retard10}
 \end{center}
 \end{figure}

   \paragraph{Commentaires sur la figure \ref{retard10}. } 
 Sur l'exemple de la figure \ref{retard5}, lorsque la trajectoire longe un segment horizontal elle traverse toujours l'axe $y=0$. Ce n'est pas  le cas lorsque le point de sortie est déterminé par :
 $$ \int_x^{S(x)}f(s)ds = 0$$
 comme le montre l'exemple de la figure \ref{retard10}.
 Sur cet exemple on a pris $\eps = 0.01$,  $f(x) = 0.5\cos(x) +0.1$ qui s'annule une première fois (pour $x \geq 0$) en $\theta_1 = \arccos (-0.2) $ et on a pris $\rho = -0.6$ de façon à ce que le graphe (en rouge)  de $x \mapsto \int_{\theta_1}^xf(s) ds$ recoupe l'axe avant de rencontrer $y = ±2\rho$. On voit que la solution qui était entrée "par le haut" au point $a$ dans le halo de $y=0$ ressort {\em par le haut} au point $b_1$.
 
\begin{figure}[ht]
  \begin{center}
 \includegraphics[width=1\textwidth]{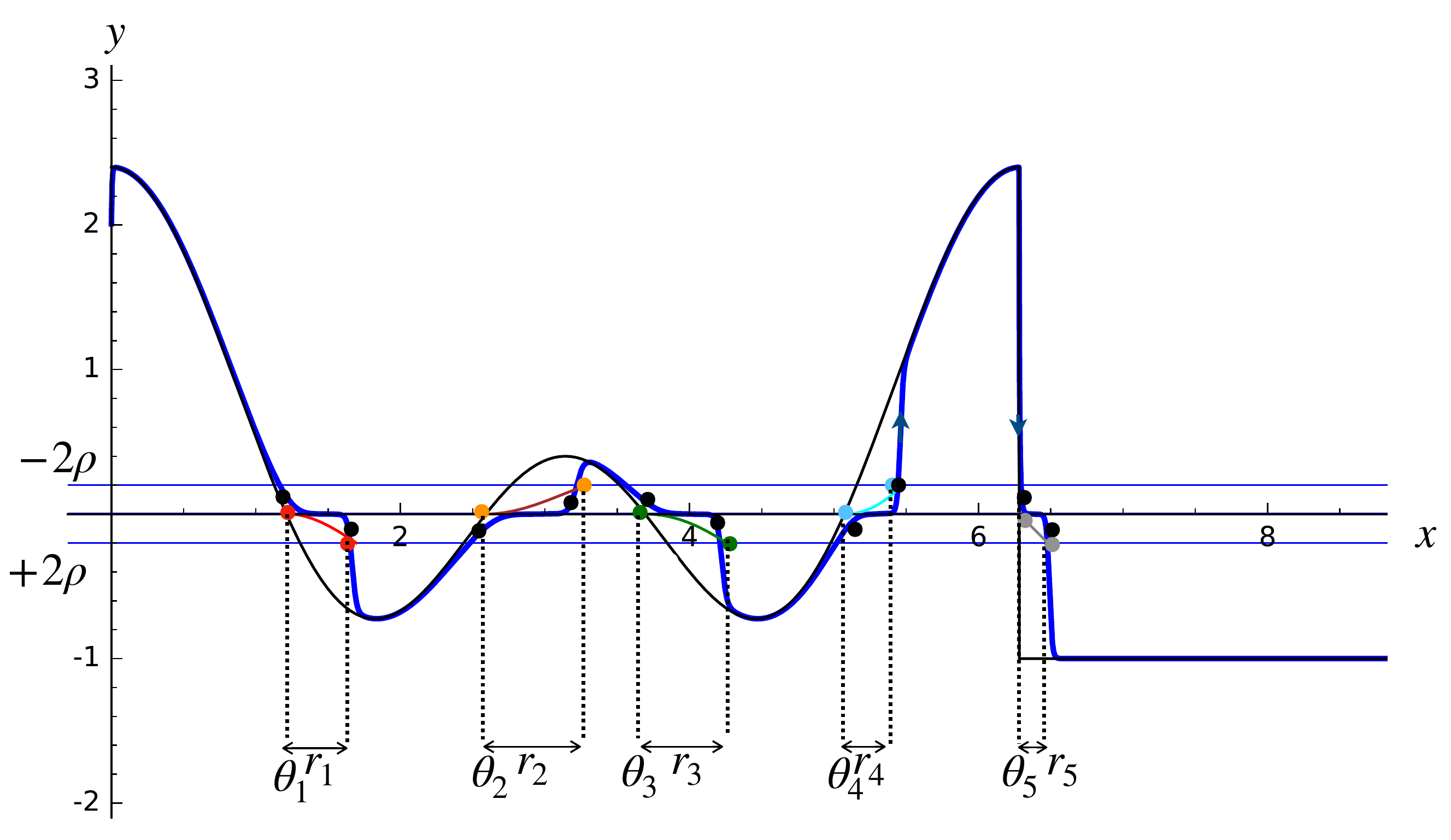}
 \caption{Le "retard à la sortie" diminue quand $\rho$ augmente. Explications dans le texte} \label{retard5bis}
 \end{center}
 \end{figure}

 \paragraph{Commentaires sur la figure \ref{retard5bis}. } 
  Pour cette simulation on a repris le modèle de la figure \ref{retard5} sauf  $\rho = - 0.1$ à la place de $\rho = -0.4$. On voit que la solution est qualitativement assez différente : elle présente maintenant 5 segments horizontaux, au lieu de 3 dans le cas $\rho = -0.4$. Le lecteur se convaincra aisément que, en toute généralité, la somme des longueurs des segments horizontaux tend vers $0$ lorsque $\rho$ tend vers $0$.
  
\begin{figure}[ht]
  \begin{center}
 \includegraphics[width=1\textwidth]{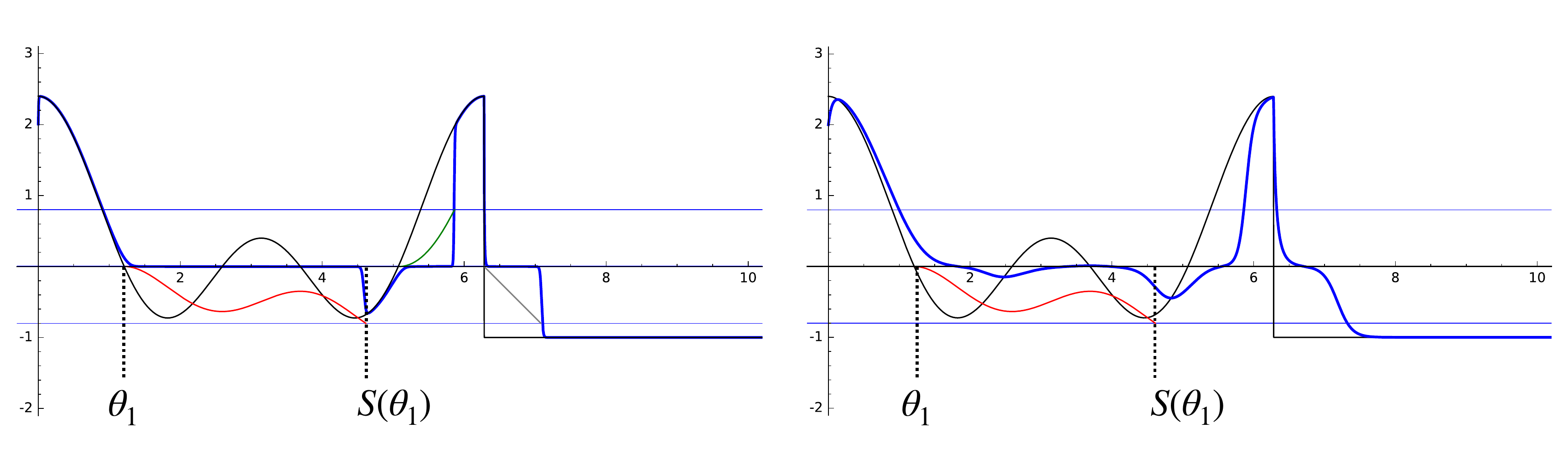}
 \caption{A gauche $\eps = 0.01$, à droite $\eps = 0.1$} \label{retard6bis}
 \end{center}
 \end{figure}

 \paragraph{Commentaires sur la figure \ref{retard6bis}.} 
  
Le théorème d'approximation est énoncé avec $\eps$ infiniment petit. Pour les simulations de la figure \ref{retard5}  j'ai pris $\eps = 0.01$. Qu'en est-il pour un ''epsilon plus gros'' ? Sur cette figure je compare la simulation obtenue pour $\eps = 0.01$ avec ce que donne la simulation du même système pour $\eps = 0.1$. Les successions de segments horizontaux, verticaux sont moins nettes mais les prédictions des points de sortie restent valables.

\newpage

\section{Démonstration du théorème d'approximation.}
Cette section est consacrée à la démonstration du théorème \ref{theoreme}. Elle consiste en l'examen de portraits de phase et le suivi de trajectoires typiques. Nous laisserons au lecteur le soin de se convaincre que tous les cas possibles ont bien été envisagés.

\subsection{En dehors du halo de l'axe horizontal $y = 0$.}
Sur la figure à droite on a représenté le graphe de $f$ (en noir) dans une zone où il ne rencontre pas l'axe $y = 0$. On considère deux conditions initiales, $(x_1,y_1)$ en dessous du graphe, non infiniment proche du graphe de $f(x)$ ni de l'axe $y = 0$ et $(x_2,y_2)$ au dessus du graphe.\\[4pt]
\begin{minipage}{0.49 \textwidth}
 Les courbes vertes sont les graphes respectifs de $f(x)+\alpha$ et $f(x)-\alpha$ où $\alpha$ est choisi strictement positif \textbf{non infiniment petit}. Les segments $[a,b]$ et $[a,c]$ sont portés respectivement par les droites de pente $-1/\alpha$ et $+ 1/\alpha$ issues du point $(x_1,y_1)$ jusqu'à leur rencontre avec le graphe de $f(x)-\alpha$. On définit ainsi un domaine $D$ délimité par les 
\end{minipage} $\quad$
\begin{minipage}{0.5 \textwidth}
\includegraphics[width=1\textwidth]{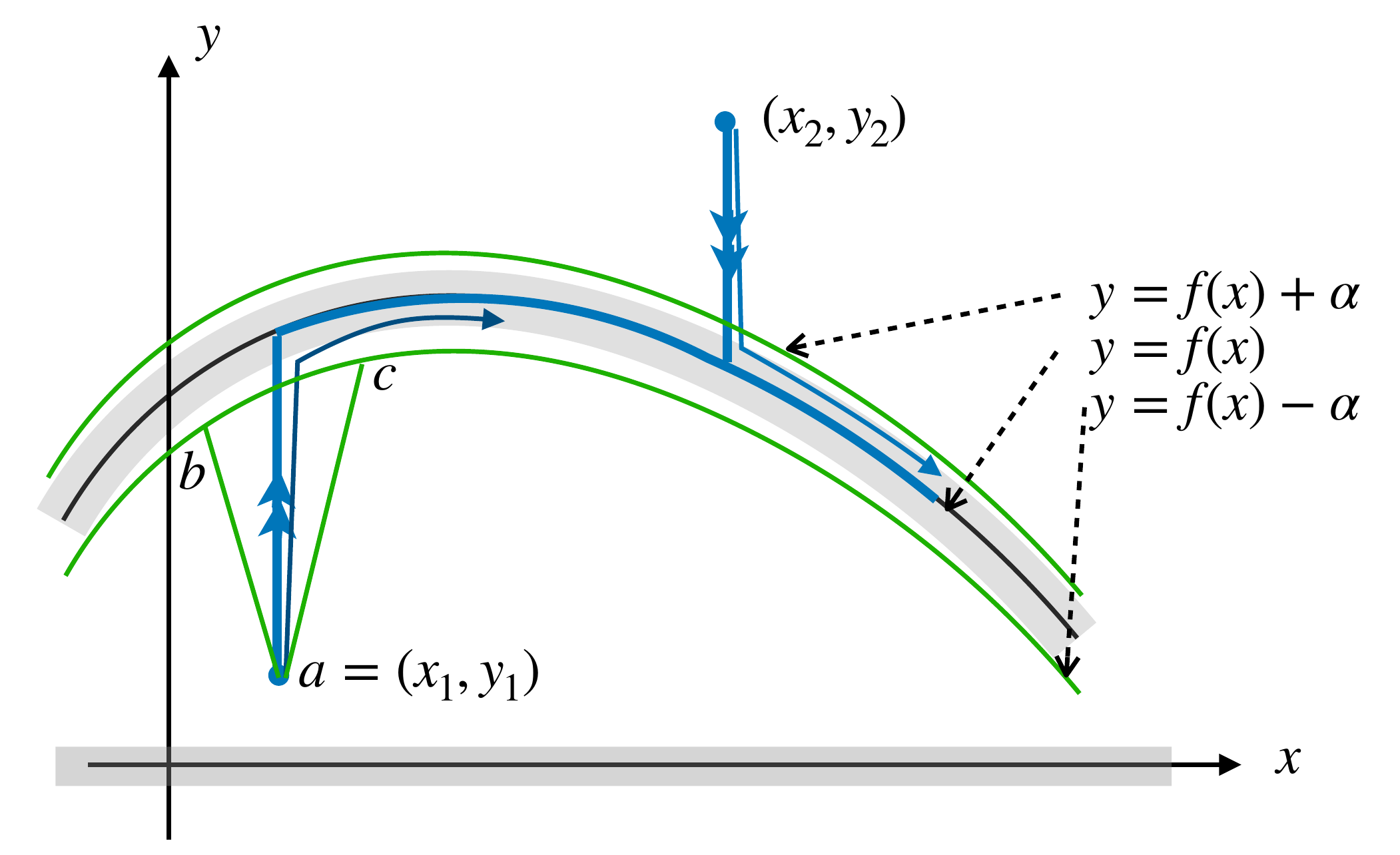}
\end{minipage} 
segments $[a,b]$, $[a,c]$ et l'arc $(a,b)$ porté par le graphe de $f(x)-\alpha$. Considérons le système $S_m$ :
 \beq \label{Sm1}
 S_m\quad  \quad \left\{
\begin{array}{lcl}
\displaystyle \frac{dx}{dt}& =& 1\\[6pt]
\displaystyle  \frac{dy}{dt} &=& \displaystyle \frac{1}{\eps}\sqrt{m^2+y^2} \Big( f(x) -y\Big)
 \end{array} 
 \right.
\feq
Comme le domaine $D$ ne rencontre pas de zone grisée, ni le facteur $\sqrt{m^2+y^2} $ ni le facteur $(f(x)-y)$ ne sont infiniment petit et donc le second membre de la seconde équation de \eqref{Sm1} est infiniment grand positif et par suite le long des segments $[a,b]$ et $[a,c]$ le champ pointe strictement à l'intérieur du domaine $D$ et, comme il ne s'annule pas, la trajectoire issue de $a = (x_1,y_1)$ ressort de $D$, au bout d'un temps nécessairement infiniment
petit, en un point de l'arc $(a,b)$. Comme ceci vaut pour tout $\alpha \gnsim 0$ il est clair que la trajectoire finira par être infiniment proche du graphe de $f(x)$, en revanche il est moins clair qu'il existe un $t_1$ \textbf{infiniment petit} pour lequel $(x(t_1),y(t_1))$ est infiniment proche du graphe de $f(x)$ car lorsqu'on se rapproche du graphe la vitesse de $y$ diminue.
Ce point est acquis par un argument typiquement nonstandard détaillé à l'annexe \ref{halocourbelente}
On se convainc aisément que, une fois que la solution  $(x(t),y(t))$ est entrée dans le halo du graphe de $f(x)$, la trajectoire doit y rester tant que $x(t)$ reste non infiniment proche d'une valeur où $f(x)$ change de signe. 

\subsection{Dans le halo de l'axe horizontal $y = 0$ : la loupe exponentielle.}
\begin{definition}
On appelle {\em loupe exponentielle} le changement de variable :
$$z = [y]^{\eps}\quad  \stackrel{df}{=} \quad \mathrm{sgn}(y)|y|^{\eps}$$
\end{definition}
Ce changement de variable a été introduit dans \cite{BEN81} et \cite{BCDD81} .\\[4pt]

\begin{figure}
  \begin{center}
 \includegraphics[width=1\textwidth]{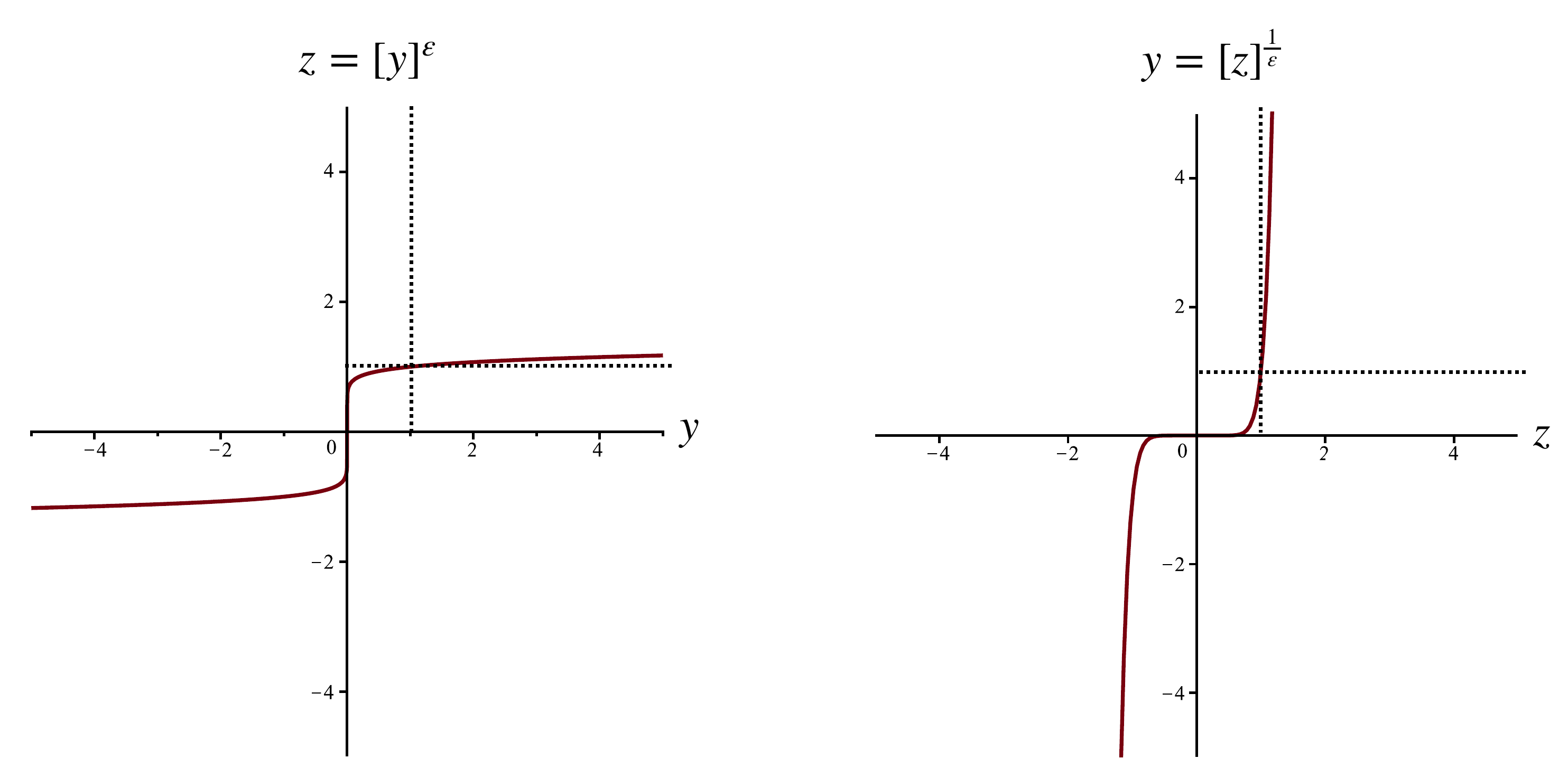}
 \caption{ La fonction $y \mapsto [y]^\eps$ est strictement croissante de $\Rmat$ dans $\Rmat$ ; sa fonction réciproque est $y = [z]^{\frac{1}{\eps}}$ et l'on a $\left([y]^{\eps} \right)' = \eps\frac{[y]^{\eps}}{y}$.
Lorsque $\eps$ est {infiniment petit}  ce changement de variable étale sur $]-1,+1[$ le {halo} de $0$ et concentre autour de $+1$ (resp. $-1$) les {limités} positifs (resp. négatifs). Ici $\eps = 0.1$. } \label{loupe}
 \end{center}
 \end{figure}
 
\noindent Si l'on fait le changement de variable $z = [y]^{\eps}$ dans le système $S_m$. 0n obtient :
\begin{equation}
\begin{array}{l}
\displaystyle \frac{dz}{dt} =\displaystyle  \eps \frac{[y]^{\eps}}{y}\frac{dy}{dt} \\[8pt]
\displaystyle\frac{dz}{dt} =\displaystyle  \eps \frac{[y]^{\eps}}{y}\left(\frac{1}{\eps}\right) \sqrt{m^2+y^2}(f(x)-y)=\frac{[y]^{\eps}}{y} \sqrt{m^2+y^2}(f(x)-y)\\[8pt]
\displaystyle\frac{dz}{dt} =\displaystyle \frac{[y]^{\eps}}{y} \sqrt{m^2+y^2}(f(x)-y) = z\,\mathrm{sgn}(y) \sqrt{1+\frac{m^2}{y^2}}(f(x)-y)\\[8pt]
\displaystyle\frac{dz}{dt} =\displaystyle |z| \sqrt{1+\frac{m^2}{y^2}}(f(x)-[z]^{\frac{1}{\eps}})
\end{array}
\feq
\textbf{Pour $ m <1$ on écrit  $m$ sous la forme :}
$$m = \exp\left(\frac{\rho}{\eps}\right)\quad \quad \rho < 0$$
ce qui donne :
$$\frac{m^2}{y^2} = \exp \left(2\frac{\rho -\ln(|z|)}{\eps}\right)$$
donc, dans les variables $(x,z)$, on étudie le système :  
\beq \label{eqlocal3}
S_m \quad \quad \quad  \left\{
\begin{array}{lcl}
\displaystyle \frac{dx}{dt}& =&1 \\[6pt]
\displaystyle \frac{dz}{dt} & =&\displaystyle |z|\sqrt{  1+\exp \left(2 \frac{\rho-\ln(|z|)}{\eps} \right)    } \left( f(x)-[z]^{\frac{1}{\eps}}\right)
\end{array} 
\right.
\feq
prolongé par continuité pour $z = 0$.

\bitbul
\item Lorsque $z \gnsim 1 $ (resp $z \lnsim -1 $) on a :
	\bito
	\item $|z| > 1 $ 
	\item $\sqrt{\cdots} > 1$ 
	\item $(f(x) - [z]^{\frac{1}{\eps}}) = -\infty$ (resp $+\infty$) parce que $[z]^{\frac{1}{\eps}}$ est {infiniment grand} positif (resp négatif)
	\fit
et donc $\displaystyle \frac{dz}{dt} = -\infty$ (resp $+\infty$)

 \item La quantité $\exp \left(2 \frac{\rho-\ln(|z|)}{\eps} \right) $ est {infiniment grande} lorsque $\rho -\ln(|z|) \gnsim 0$ et donc 
 $$ |z|\lnsim \emat^{\rho}\Longrightarrow   \frac{dz}{dt} = \mathrm{sgn}(f(x))\cdot(+\infty)$$

 \item Dans les bandes $\Rmat \times \eset{z:\;-1 \lnsim z \lnsim - \emat^{\rho}}$ et  $\Rmat \times \eset{z :\; \emat^{\rho} \lnsim z \lnsim 1} $ on a $\rho+\ln(|z|) \gnsim 0$ et donc $\exp \left(-2 \frac{\rho+\ln(|z|)}{\eps} \right)  \sim 0$. Donc dans ces deux bandes les trajectoires sont infiniment proches des trajectoires  du champ {standard} : 
 \beq \label{eqlocal4}
\stackrel{\sim}{S_m} \quad \quad \quad \quad \left\{
\begin{array}{lcl}
\displaystyle \frac{dx}{dt}& =&1 \\[6pt]
\displaystyle \frac{dz}{dt} & =& |z]f(x)
\end{array} 
\right.
\feq
 \fit
 A partir de ces informations il est possible de faire une esquisse du portrait de phase de $S_m$ dans les variables $(x,z)$, ce que nous entreprenons dans la sous-section suivante.

\subsubsection{Commentaires sur la figure \ref{schema2} }
\begin{figure}[ht]
  \begin{center}
 \includegraphics[width=0.9\textwidth]{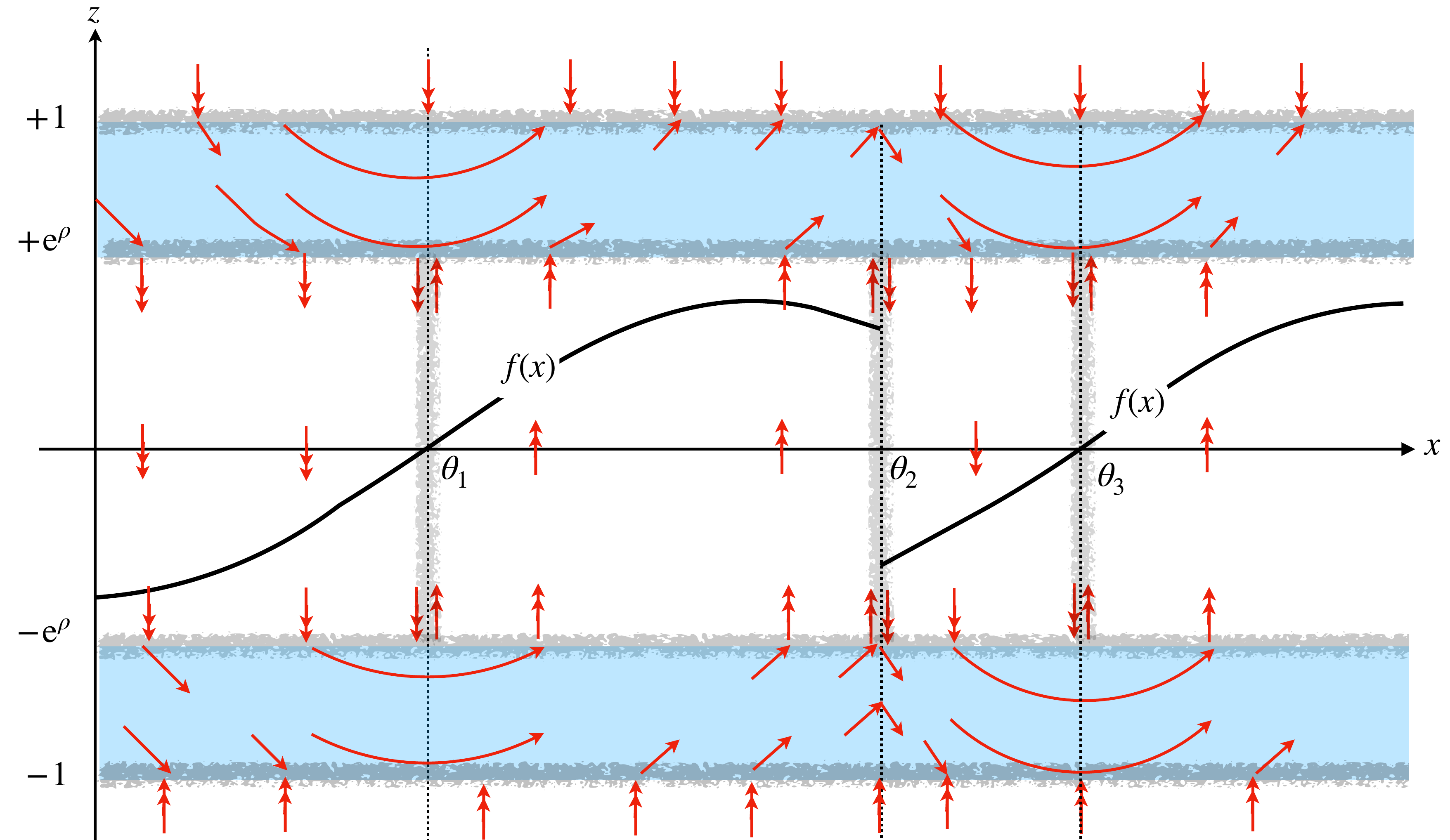}
 \caption{Le système $S_m$ sous la loupe exponentielle. Attention le graphe de $f$ est dans les variables originales $(x,y)$} \label{schema2}
 \end{center}
 \end{figure}
\ben
\item La figure laisse apparaître une partition en zones délimitées par des horizontales et des verticales.
 Les verticales pointillées sont les droite $x = \theta_n \;n\in J$ où la fontion $f$ change de signe.

\item Les horizontales épaisses en grisé représentent les {halos} des droites $z = -1$, $z = -\emat^{\rho}$, $z = \emat^{\rho}$ et $z = 1$ ;

\item Les variables sont $(x,z)$ mais le graphe de $f$ est représenté dans les variables $(x,y)$ ; il n'est là que pour matérialiser les valeurs de $x$ où la fonction $f$ change de signe ; sur cet exemple  il y a trois changements de signe (pour $x = \theta_1,\,\theta_2,\,\theta_3$)  ; le changement de signe du milieu correspond  à une discontinuité, les deux autres correspondent  à des zéros de $f$. De même que ci-dessus, on a symbolisé en gris les {halos} des segments verticaux $x = \theta_n $ compris entre $-1$ et $+1$.

\item Dans les zones non grises,  la seconde composante du champ $S_m$  est de valeur absolue {infiniment grande} positive ou négative. Ceci est symbolisé par des doubles flèches verticales orientées ; au dessus de $+1$ elles sont toujours dirigées vers le bas, au dessous de $-1$, dirigées vers le haut et dans la bande $ -\emat^{\rho}< z <  \emat^{\rho}$ selon le signe de $f(x)$.

\item Dans les zones bleues on a symbolisé (en rouge)  les variations des trajectoires du champ  $\stackrel{\sim}{S_m}$ au voisinage des changement de zone.

\item Dans les zones grises la trajectoire évolue vers la droite à vitesse unité. En particulier les zones grises verticales sont traversées en une durée infinitésimale.
\fen
 \clearpage
 
 \subsubsection{Commentaires sur la figure \ref{schema2-2} : \\Suivi de trajectoires dans les variables $(x,z)$}\label{espacexz}
\begin{figure}[ht]
  \begin{center}
 \includegraphics[width=0.93\textwidth]{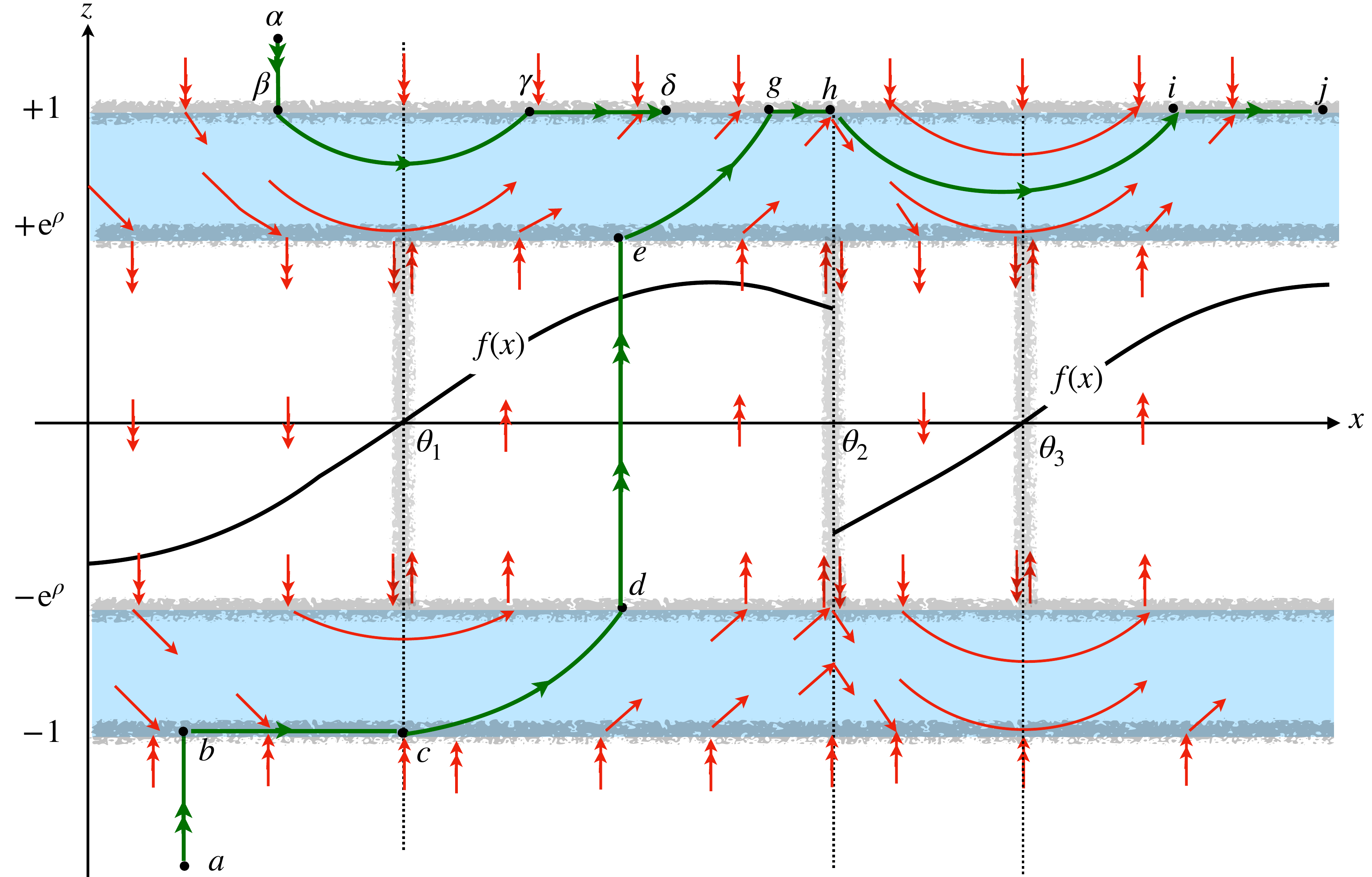}
 \caption{Explication paragraphe \ref{espacexz}.} \label{schema2-2}
 \end{center}
 \end{figure}
On considère une condition initiale $a = (x_0,z_0)$ avec $x_0 >0$ et $z_0 \lnsim -1$ non infiniment grande, comme indiqué sur la figure ;  le point $(x_0,-1)$ est attractif (i.e. au dessus le champ est dirigé vers le bas, vers le haut en dessous) ; on note $(x(t),z(t))$ la trajectoire correspondante. On a évidemment $x(t) = x_0+t$.
\paragraph{Trajectoire $a,b,c,d,e,g,h,i,j$:}
\ben
\item \texttt{Segment} $\arc{a,b}$. Le champ est infiniment grand dirigé vers le haut. La trajectoire correspondante est quasiment verticale et il 
 existe un $t _1\sim 0$ tel que $z(t_1)\sim-1$ (voir appendice \ref{demiLR}).
 
\item \texttt{Segment} $\arc{b,c}$. Puisque $z = -1$ est attractif  jusqu'au point $c$ la trajectoire {longe} la droite $z = -1$ à la vitesse $+1$ jusqu'au point $c$ (voir appendice \ref{demiLR}) qu'elle atteint au temps $t_2 \sim \theta_1-x_0$. 

\item \texttt{Segment} $\arc{c,d}$. Il existe un instant $t_3\sim t_2 $ tel que $(x(t_3),z(t_3))\sim c$ soit dans la zone bleue où $S_m$ est {infiniment proche} de $\stackrel{\sim}{S_m}$. Donc le segment $\arc{c,d}$ est {infiniment proche}  de la trajectoire $(\tilde{x}(t), \tilde{z}(t)$ de $\stackrel{\sim}{S_m}$ issue de $c$ tant que $\tilde{z}(t) \leq  - \emat^{\rho}$ (voir appendice \ref{changementdesigne} pour une preuve plus formelle).
On a :
\beq
\begin{array}{lcl}
\tilde{x}(t )&=& \theta_1 +t\\[6pt]
\tilde{z}(t) &=& \displaystyle -\emat^{- \int _{\theta_1} ^t f(s) ds} 
\end{array}
\feq
et donc la valeur $-\emat^{\rho}$ est atteinte au temps $\tau_4$ tel que :
\beq 
\displaystyle -\emat^{- \int _{\theta_1} ^{\tau_4}f(s) ds} = -\emat^{\rho}
\feq 
soit :
\beq \label{tau1}
\displaystyle \int _{\theta_1} ^{\tau_4}f(s) ds = -\rho
\feq 
Il existe un $t_4 \sim \tau_4$ tel que $ z(t_4) \sim -\emat^{\rho}  $

\item \texttt{Segment} $\arc{d,e}$. A partir de $t_4$ la trajectoire est quasi-verticale ascendante. Il existe $t_5 \sim t_4$ tel que $z(t_5) \sim \emat^{\rho} $. Au point $e$ le champ est traversant. 

\item \texttt{Segment} $\arc{e,g}$. A partir de $t_5$ nous sommes dans la même situation qu'au temps $t_3$ (voir 3. ci-dessus). La trajectoire de $\stackrel{\sim}{S_m}$ atteint donc la valeur $1$ au temps $\tau_6$ tel que :
\beq 
\displaystyle \emat^{\rho} \emat^{ \int _{\tau_5} ^{\tau_6}f(s) ds} = 1
\feq 
soit :
\beq \label{tau2}
\displaystyle \int _{\tau_5} ^{\tau_6}f(s) ds = -\rho
\feq 

\item En réunissant \eqref{tau1} et \eqref{tau2} il vient que $t_6 \sim \tau_6$ où $\tau_6$ est défini comme le premier instant tel que :
\beq \label{condition1}
\displaystyle \int_{\theta_1} ^{\tau_6}f(s) ds = - 2\rho
\feq 

\item \texttt{Segment} $\arc{g,h}$. Pour les mêmes raisons que 2. ci-dessus le segment $\arc{f,g}$ est dans le {halo} de $z = +1$ tant que $t \lnsim \theta_2$.

\item \texttt{Segment} $\arc{h,i}$. Pour les mêmes raisons que 3. ci-dessus la trajectoire   de  $S_m$  reste dans la zone bleue jusqu'au temps $t_7 \sim \tau_7$ défini par :
\beq \label{condition2}
\displaystyle \int_{\theta_2} ^{\tau_7}f(s) ds =  0
\feq 
A la différence du segment $\arc{c,d}$ qui traversait la zone bleue entre $z = - 1$ et $z = -\emat^ {\rho} $ ici le segment $\arc{h,i}$, issu d'un point tel que $z = 1$, remonte vers un point tel que $z = 1$ sans traverser la zone bleue.

\item \texttt{Segment} $\arc{i,j}$. Comme pour $\arc{b,c}$ on reste dans le {halo} de $z = 1$.
\fen
\paragraph{Trajectoire $\alpha,\beta,\gamma,\delta$:}
Alors que la trajectoire précédente traverse la bande $-1< z<  +1$ cette trajectoire pénètre dans la zone bleue mais ne traverse pas la bande ni l'axe $y = 0$.

\begin{remarque}
Il n'y a pas de différence qualitative entre le point $g$ qui correspond à une discontinuité de $f$ et le point $c$ qui correspond à un zéro de $f$ à part le caractère non $C^1$ de la trajectoire lors du changement de signe.
\end{remarque}

 \subsection{Suivi d'une trajectoire : \\Utilisation simultanée des variables $x,y$ et $x,z$.}
\begin{figure}[ht]
  \begin{center}
 \includegraphics[width=1\textwidth]{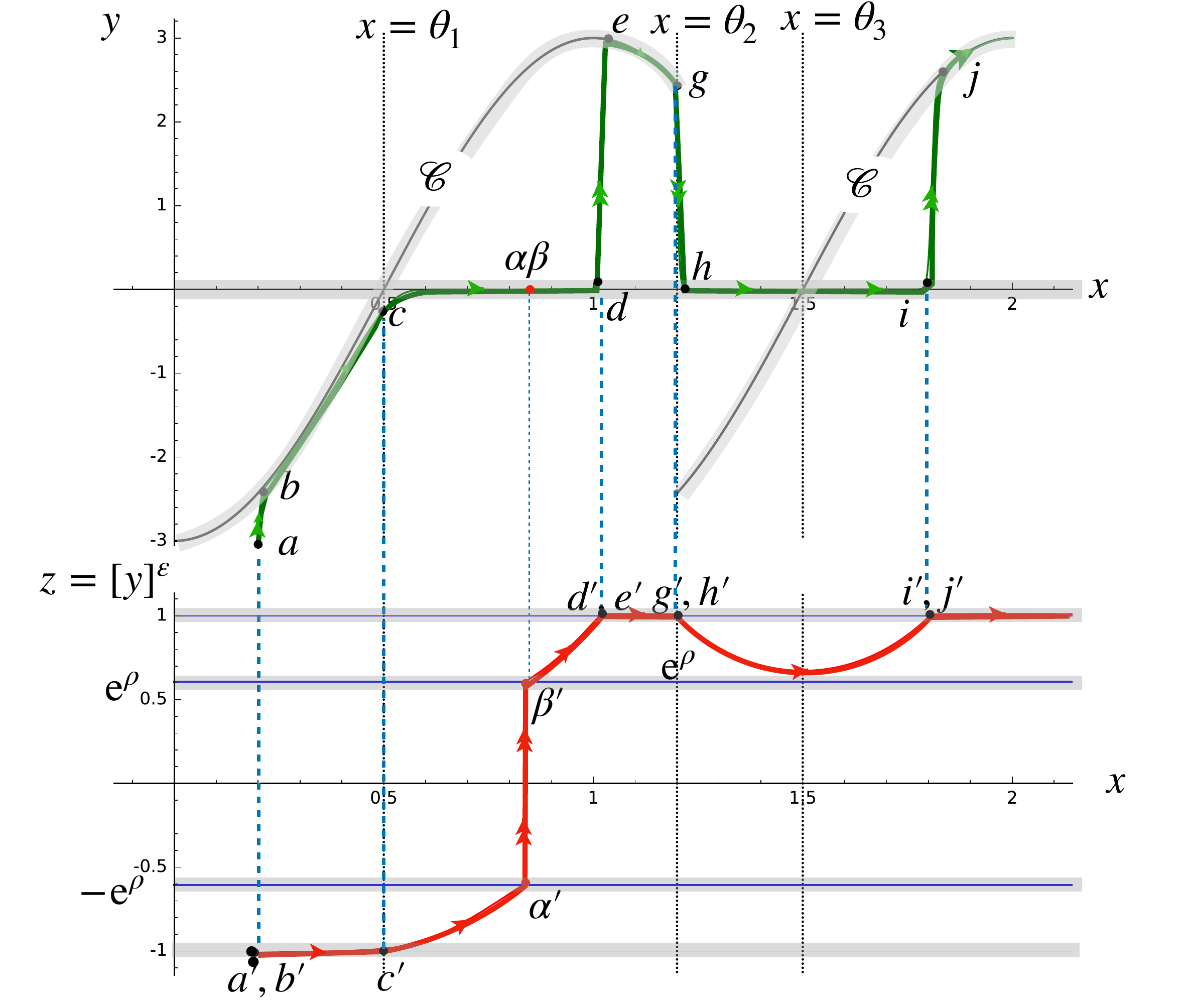}
 \caption{Au dessus le système $S_m$ dans les variables $x,y$ et au dessous l'agrandissement du halo de l'axe $y=0$ sous la loupe exponentielle. sous la loupe exponentielle.} \label{loupesimu1}
 \end{center}
 \end{figure}
 
 Sur la figure \ref{loupesimu1} sont superposés les espaces des trajectoires de $S_m$, (avec $m= \emat^{\frac{\rho}{\eps}}$) dans les variables $(x,y)$ et les variables $(x,z)$, les valeurs de $x$ se correspondant. 
 \bit
 \item Dans les variables $(x,y)$ on a représenté le graphe (courbe $\mathcal{C}$) de la fonction $f$ utilisée dans la section précédente,  qui s'annule deux fois en $x =\theta_1$ et $x =\theta_3$  et possède une discontinuité où elle change de signe en $x =\theta_2$. les halos de $\mathcal{C}$ et de l'axe $y = 0$ sont symbolisés en grisé.
 
  \item Dans les variables $(x,z)$ on a symbolisé en grisé les halos des horizontales $z = ±1$ et  $z = ± \emat^ {\rho}$.
 \fit 
 Nous allons suivre la trajectoire issue du point $a$.
 \bitbul
 \item \texttt{Segment $\arc{a,b}$}. La courbe lente est attractive. On longe un segment vertical puis on pénètre dans le halo de $\mathcal{C}$ en un point $b$ d'abscisse infiniment proche de celle de $a$. Dans les variables $(x,z)$ les deux points correspondants, $a'$ et $b'$ sont dans le halo de $z = - 1$ et infiniment proches l'un de l'autre.

 \item \texttt{Segment $\arc{b,c}$}. La trajectoire longe la courbe lente attractive tant que $f(x) \lnsim 0$ jusqu'au point $c = (x_1,y_1)$ où elle pénètre dans le halo de $y = 0$. Dans l'espace $(x,z)$ le segment $\arc{b',c'}$ reste dans le halo de $z = -1$ jusqu'au point $c'$ infiniment proche de $(\theta_1,-1)$.
\fit 
Dans l'espace $(x,y)$, à partir du point $c$ on a deux possibilités : rester dans le halo de $\mathcal{C}$ ou bien dans celui de $y = 0$ ; la réponse est donnée dans l'espace $(x,z) $. On suit maintenant la trajectoire dans l'espace $(x,z)$ à partir du point $c'$ jusqu'au point $d'$.
\bitbul
 \item \texttt{Segment $\arc{c',\alpha'}$}.  On suit maintenant la trajectoire dans l'espace $(x,z)$ à partir du point $c'$  jusqu'à l'entrée dans le halo de la droite $z = -\emat^{\rho}$ au point $\alpha'$. Dans l'espace $(x,y)$, puisque $|z| \lnsim 1$, le segment $\arc{c,\alpha}$ reste dans le halo de $y = 0$.

 \item  \texttt{Segment $\arc{\alpha',\beta'}$}. Dans l'espace $(x,z)$ le segment $\arc{\alpha',\beta'}$ fait sauter en une durée infiniment petite de $z =-\emat^{\rho}$ à $z =+\emat^{\rho}$  (paragraphe \ref{espacexz}, point 3.) Dans l'espace $(x,y)$ le segment $\arc{\alpha,\beta}$ est infiniment court et traverse l'axe $y = 0$.
 
  \item \texttt{Segment $\arc{\beta',d'}$}. On suit  la trajectoire partir du point $\beta'$  jusqu'à la rencontre avec le halo de  la droite $z = +1$ au point $d'$. Pénétrer dans le halo de $z = 1$ c'est quitter le halo de $y = O$, donc, dans l'espace $(x,y)$ le point $d = (x_2,y_2) $ correspondant est le point de sortie du halo de $y = 0$. La valeur de $x_2$ est donnée par :
  $$\int_{\theta_1}^{x_2}f(s)ds   =- 2\rho$$   
  \fit
    Puisqu'on est sorti du halo de $y = 0$  on peut retourner dans l'espace des $(x,y)$.
  \bitbul
  
 \item \texttt{Segment $\arc{d,e,g,h}$}. On saute en $e$ où l'on pénètre dans le halo de la courbe $\mathcal{C}$ qu'on longe jusqu'au point $g$ où l'on saute vers le halo de $y =0$ où l'on pénètre au point $h$. Le segment correspondant dans l'espace $(x,z)$, $\arc{d',e,',g',h'}$ longe $z = 1$.
 \fit
 On retourne dans l'espace $(x,z)$.
 \bitbul
 
  \item \texttt{Segment $\arc{h',i'}$}. Dans l'espace des $(y,z)$ on suit la trajectoire qui ne traverse pas la bande $ \emat^{\rho} < z < 1$ mais rejoint à nouveau le halo de $z =1$ au point $i''$. Dans l'espace $(x,y)$, le point correspondant $h =(x_3,y_3)$  est le point de sortie du halo de $y = 0$. La valeur de $x_3$ est donnée par :
  $$\int_{\theta_1}^{x_3}f(s)ds   = 0$$
 
 \fit
 On retourne dans l'espace  $(x,y)$.
 \bitbul
  \item \texttt{Segment $\arc{i,j}$}. Dans l'espace des $(x,y)$ on saute de $i$ au point $j$ dans le halo de $\mathcal{C}$.
 \fit

 \subsubsection{Condition initiale dans le halo de l'axe $y=0$}\label{ci0}
  \begin{figure}[ht]
  \begin{center}
 \includegraphics[width=1\textwidth]{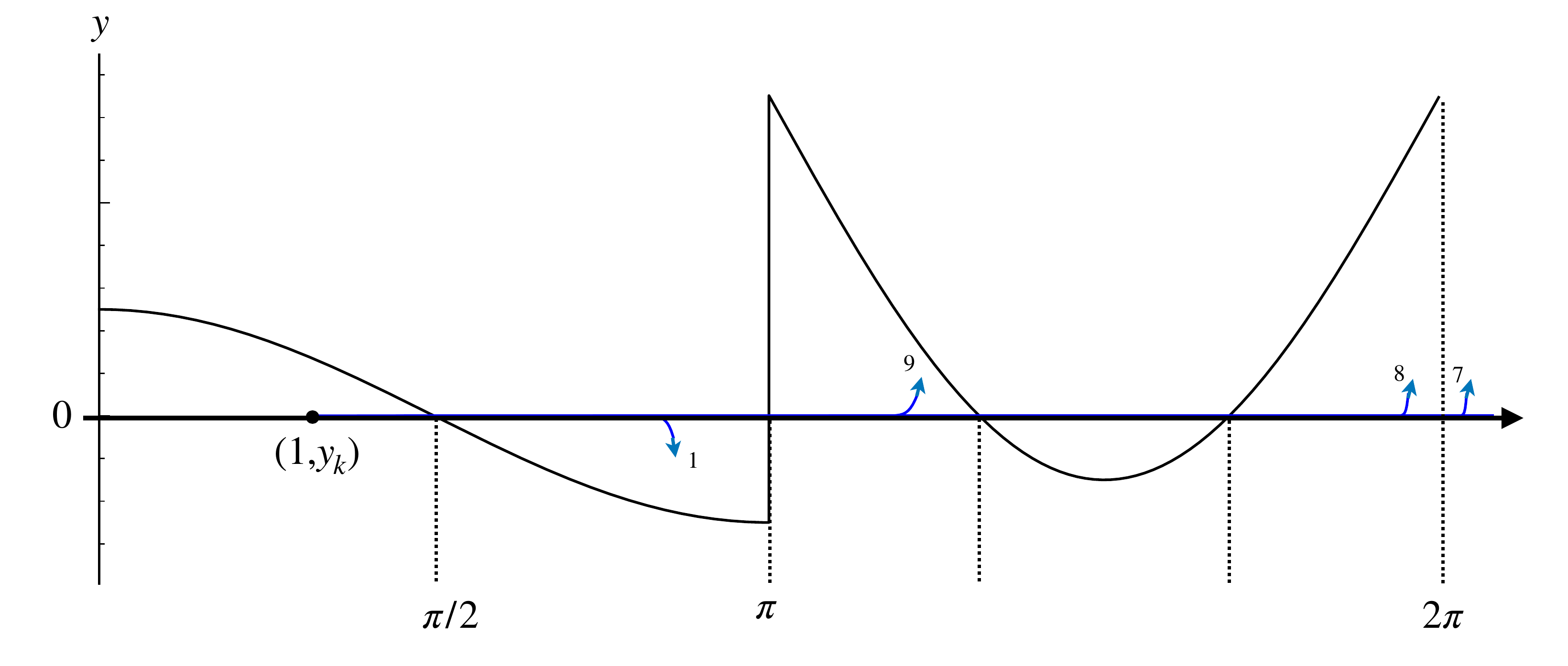}
 \caption{Simulation de 9 trajectoires issues de 9 points infiniment proches. Explications paragraphe \ref{ci0} } \label{cisim01}
 \end{center}
 \end{figure}
  \begin{figure}[ht]
  \begin{center}
 \includegraphics[width=1\textwidth]{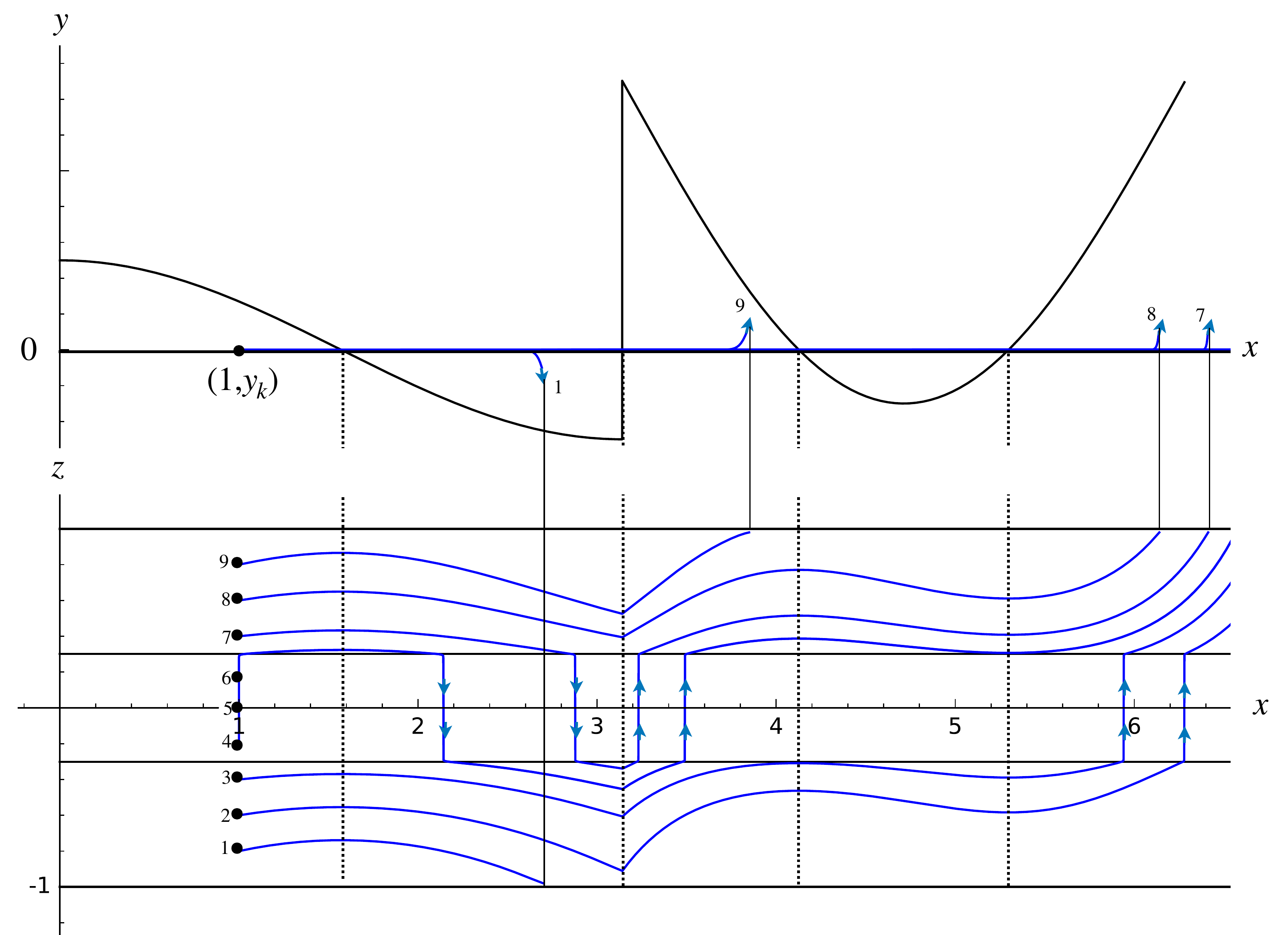}
 \caption{Agrandissement sous la loupe exponentielle des trajectoires de la simulation \ref{cisim01}. Explications paragraphe \ref{ci0} } \label{cisim02}
 \end{center}
 \end{figure}
 Soit une condition initiale $(x_0, y_0)$ telle que $y_0 \sim 0$ (condition exclue par le théorème d'approximation). Que pouvons nous dire de la solution issue de $(x_0, y_0))$ sans autre information sur $y_0$ ? 
 
 Sur la figure \ref{cisim01} on observe le résultat de la  simulation de $S_m$ avec la fonction $f$ donnée par :
 \beq \label{fcisim0}
\begin{array}{lcl}
\displaystyle x \in [0,\pi[ &\Longrightarrow &f(x) = 0.5\cos(x) \\[8pt]
\displaystyle x \in [\pi,2 \pi[& \Longrightarrow &f(x) =  1.5 + 1.8\sin(x)
\end{array}
\feq
 dont le graphe est tracé en noir,  pour les valeurs des paramètres :
 $$ \eps = 0.01 \quad \quad \quad m = \emat^ {\rho/\eps} \quad \rho = -1.2$$
 et à partir des $9$ conditions initiales :
 $$ (x_0 = 1, y_k)\quad \quad y_k = (k-5)\times 0.2 \times m \quad k = 1,2,\cdots, 9$$
 sur une durée de 5.5 unités. Les conditions initiales sont très rapprochées ; avec les valeurs des paramètres on a $|y_k|< 2\times 10^{-10}$, les 9 conditions initiales sont donc indiscernables sur la figure.
 
 La simulation montre que sur les $9$ trajectoires, seulement $4$ (correspondant à $k = 1,7,8,9$) ont quitté le halo de $0$ avant le temps $t =5.5$ en des points différents.
 
 L'application de la loupe exponentielle (dans la simulation au lieu d'afficher $y$ on affiche  $z = \emat^{\eps\log(y)}$ si $y >0$ et $z = -\emat^{\eps\log(-y)}$ si $y <0$) sépare les conditions initiales et les trajectoires associées. C'est ce qui est montré sur la figure \ref{cisim02} où l'on voit que le point de sortie du halo de l'axe $y = 0$ dépend de façon essentielle de la distance (infiniment petite) de $y_k$ à $0$. C'est la raison pour laquelle, dans le théorème d'approximation on exclu comme condition initiales possibles les points du halo de l'axe $y = 0$. En revanche, lorsqu'une trajectoire pénètre dans le halo de l'axe $y = 0$ pour une valeur proche de $x$, elle ne peut le faire qu'en passant par les valeurs $z \sim ± 1$ ce qui supprime toute incertitude sur la valeur de $S_{\rho}(x)$.

 \subsection{Feuilletage des trajectoires lorsque $f$ change de signe.}
  \begin{figure}[ht]
  \begin{center}
 \includegraphics[width=1\textwidth]{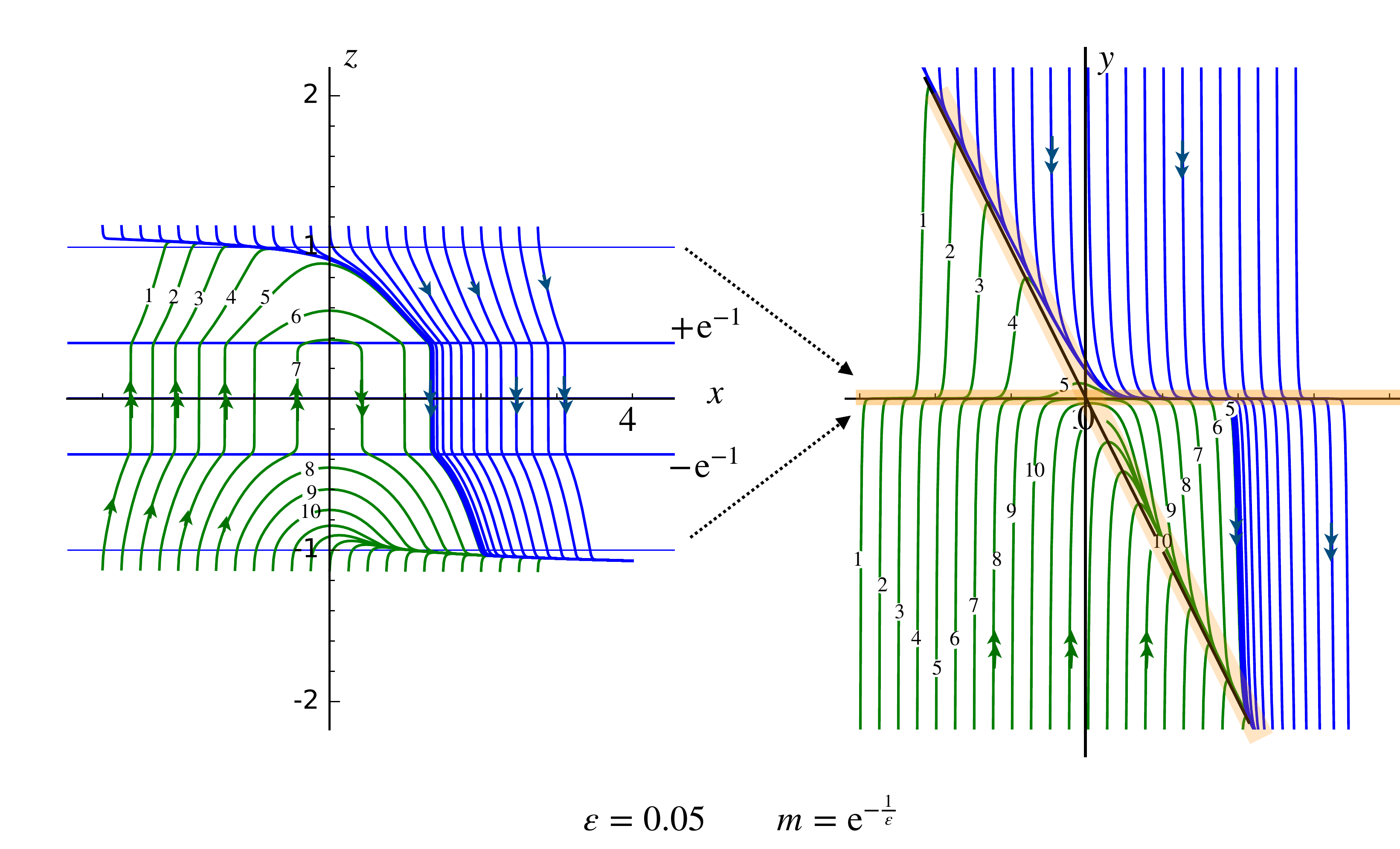}
 \caption{Simulation de $S_m$ avec $f(x) = -x$} \label{retard4}
 \end{center}
 \end{figure}
  \begin{figure}[ht]
  \begin{center}
 \includegraphics[width=1\textwidth]{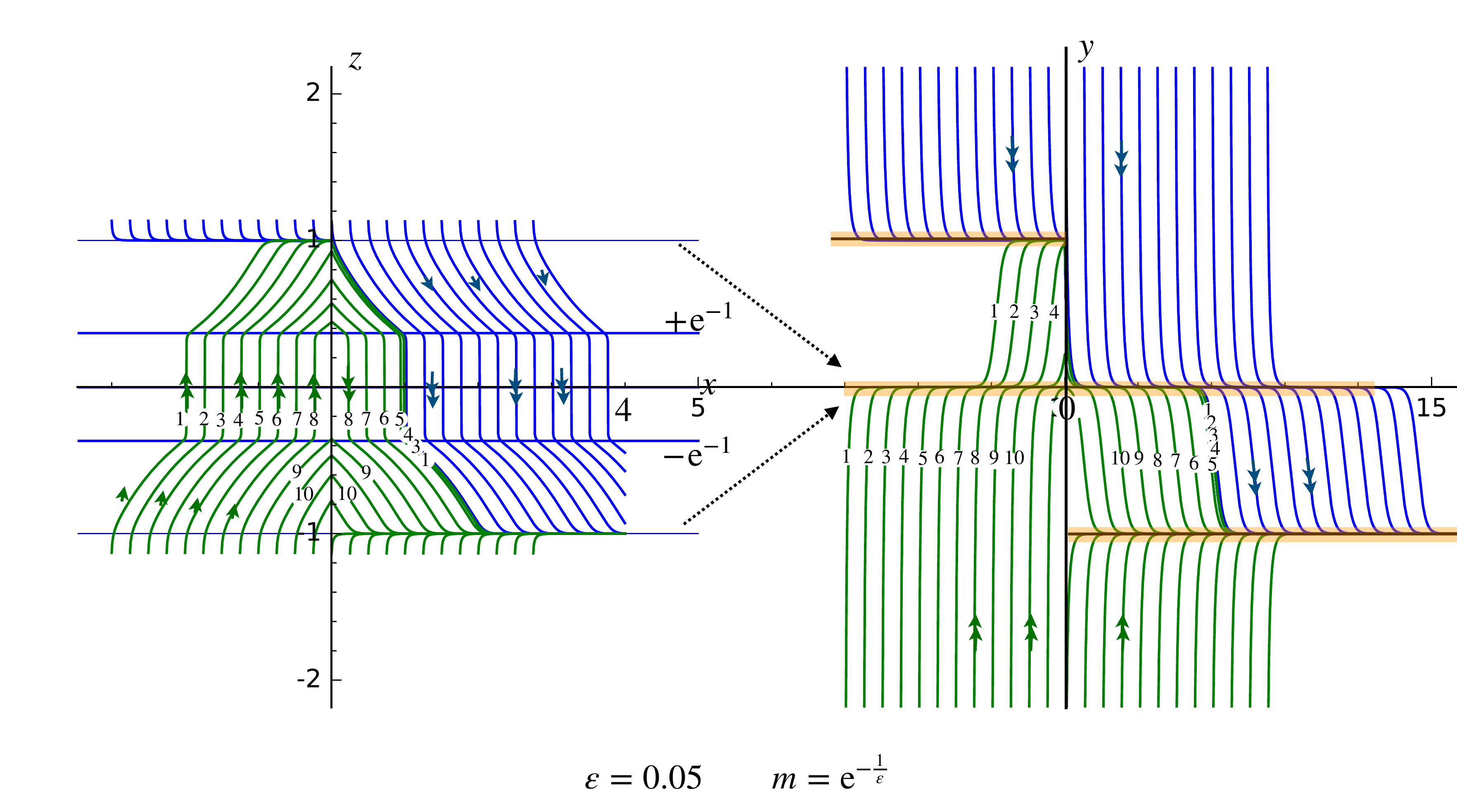}
 \caption{Simulation de $S_m$ avec $f(x) = ±1$} \label{retard4bis}
 \end{center}
 \end{figure}
 Sur les figures \ref{retard4} et  \ref{retard4bis} on a simulé les trajectoires du système $S_m$ autour de $(0,0)$ dans le cas où $f(x) = -x$ et où $f(x) = 1$ si $x$ est négatif, $-1$ sinon,  pour illustrer le passage d'un changement de signe de $f$ lorsque $f$ s'annule ou est discontinue. A droite nous avons les trajectoires dans les variables originales où les halos de $y = 0$ et $y = f(x)$ sont symbolisés par les bandes jaunes et à gauche sous la loupe exponentielle.
 
\clearpage
\newpage
 \section{Démonstration d'une conjecture de G. Katriel}
\subsection{Le problème }        
            
 On considère le système :            
 \beq \label{Kat1}
\begin{array}{lcl}
\displaystyle  \frac{dx_1}{dt} &=&r_1(\nu t)x_1+\mu(x_2-x_1)\\[8pt]
\displaystyle  \frac{dx_2}{dt} &=& r_2(\nu t) x_1+\mu(x_1-x_2)
 \end{array} 
\feq
où les fonctions $r_i$ sont standard, périodiques de période $2 \pi$ et différentiables par morceaux au sens de la section \ref{equationsdebase}. On s'intéresse au comportement de la quantité :
\beq \label{Delta}
\displaystyle \Delta(\nu,\mu) = \lim_{t \rightarrow \infty} \frac{\ln(x_1(t)x_2(t))}{t} 
\feq
lorsque :
\beq \label{HK}
\displaystyle \frac{1}{2\pi} \int_0^{2\pi} r_i( t)dt < 0 \quad \quad \chi =  \frac{1}{2\pi} \int_0^{2\pi} \max_{i = 1,2} r_i(t)dt  >0
\feq
Le système \eqref{Kat1} représente la croissance de deux populations de taille $x_1$ et $x_2$, sur deux sites $1$ et $2$, soumises à des environnements différents, $r_i(\nu t)$ en présence de migration entre les deux sites. Les hypothèses sur les $r_i$ ont pour conséquence que, en l'absence de migration ($\mu = 0$), sur chacun des sites les populations disparaissent. Appelons {\em seuil d'inflation} et notons $\mu^*(\nu)$ :
\beq \label{seuil }
\displaystyle  \mu^*(\nu) = \inf_{\mu \geq 0} \{ \mu :  \Delta(\nu,\mu) > 0 \}
\feq

Dans \cite{KAT21} il est conjecturé que :
\begin{proposition}\label{conjecturekat}  Lorsque $\nu$ est infiniment petit (la période de variation de l'environnement est très grande) le seuil de l'inflation est de l'orde de $\exp(-\frac{1}{\nu})$ (exponentiellement petit par rapport à $\nu$).
\end{proposition} 

La version classique de cette conjecture est que :

\begin{proposition}\label{conjecturekatclassique}
On considère le système \eqref{Kat1} où les fonctions $r_i$ sont périodiques de période $2 \pi$ et différentiables par morceaux au sens de la section \ref{equationsdebase}. Il existe une constante $C>0$ telle que pour $\nu$ assez grand $\mu^*(\nu) < C\exp(-\frac{1}{\nu})$.
\end{proposition} 

Je démontre cette conjecture au paragraphe suivant. Je renvoie à \cite{BLSS21, HOLTPNAS20, KLA08} pour plus d'informations sur le phénomène d'inflation mis en évidence dans \cite{HOLT02}.

\subsection{Réduction à un système $S_m$ } 
\paragraph{Passage aux variables $U = \ln(x_1x_2)$ et $V = \ln(x_1/x_2)$}$\,$\\
Posons :
$$ \xi_1 = \ln(x_1)\quad \quad \xi_2 = \ln(x_2)$$
Dans ces variables   le système \eqref{Kat1}  devient  :
  \beq \label{Kat2}
\begin{array}{lcl}
\displaystyle  \frac{d\xi_1}{dt} &=&r_1(\nu t)+\mu\left(\emat^{\xi_2-\xi_1}-1\right)\\[10pt]
\displaystyle  \frac{d\xi_2}{dt} &=& r_2(\nu t) )+\mu\left(\emat^{\xi_1-\xi_2}-1\right) 
 \end{array} 
\feq
et maintenant on pose  :
$$ U = \xi_1+\xi_2\quad \quad V = \xi_1 - \xi_2$$
Dans ces variables le système \eqref{Kat2}  devient :
 \beq \label{Kat3} 
\begin{array}{lcl}
\displaystyle  \frac{dU}{dt} &=&\displaystyle  r_1(\nu t)+r_2(\nu t) + 2\mu (\ch(V)-1)    \\[8pt]
\displaystyle  \frac{dV}{dt} &=&(r_1(\nu t)-r_2(\nu (t))-2\mu\,\sh(V)
 \end{array} 
\feq
On voit que la variable $V$ est découplée de $U$. On a :
$$\Delta(\nu,\mu) = \lim_{t \rightarrow \infty} \frac{U(t)}{t} $$
Soit  l'équation :
\beq \label{V1}
  \frac{dV}{dt} =(r_1(\nu t)-r_2(\nu (t))-2\mu\,\sh(V)
\quad \quad 
\feq
C'est une équation non autonome périodique.
\begin{proposition}
Lorsque $\mu > 0$ l'équation  \eqref{V1} possède une unique solution périodique (de période $\frac{2\pi}{\nu}$) globalement asymptotiquement stable $V_{\nu,\mu}(t)$.
\end{proposition}
\textbf{Preuve} Pour $M$ suffisamment grand on a pour $|V| = M$, $ \displaystyle \mathrm{sgn}\left(\frac{dV}{dt} \right)  =-\mathrm{sgn} (V)$ ce qui entraine que $[-M,+M]$ est invariant, l'application $V_0 \mapsto V(2\pi/\nu, 0,V_0)$ possède au moins un point fixe.
Pour un tel point fixe, l'équation aux variations le long de la solution périodique associé $V_{\nu,\mu}(t)$  est :
$$\frac{d\delta V(t)}{dt}  = -2\mu\cosh(V_{\nu,\mu} (t))\delta V(t)$$
Donc la dérivée de $V_0 \mapsto V(2\pi/\nu, 0,V_0)$ est positive strictement plus petite que $1$ ce qui prouve l'unicité et l'attractivité du point fixe.\\
$\Box$\\
Il en découle que :
$$\Delta(\nu,\mu) = \lim_{t \rightarrow +\infty} \frac{U(t)}{t} = \frac{\nu}{2\pi} \int_0^{2\pi/\nu}\displaystyle  r_1(\nu s)+r_2(\nu s) + 2\mu \left(\ch(V_{\nu,\mu}(s))-1  \right) ds $$

\paragraph{ Réduction à une période fixe.} 0n se ramène à la période $2\pi$ en posant $\overline{V} (t) = V(t/\nu)$ ce qui donne :
\beq \label{V2}
  \frac{d\overline{V}}{dt} =\frac{1}{\nu} \big(r_1(t)-r_2(t)-2\mu\,\sh(\overline{V})\big)
\quad \quad 
\feq
et : 
$$\Delta(\nu,\mu) = \frac{1}{2\pi} \int_0^{2\pi}\displaystyle  r_1( s)+r_2(s) + 2\mu \left(\ch(\overline{V}_{\nu,\mu}(s))-1  \right) ds $$
où $\overline{V}_{\nu,\mu}$ est la solution périodique de \eqref{V2}.
\paragraph{Linéarisation du sinus hyperbolique.} On fixe $\nu = \eps \sim 0$ une fois pour toutes et $\mu$ est un paramètre. On pose :
$$ W = 2\mu\sinh(\overline{V})$$
ce qui donne :
\beq \label{W1}
\frac{dW}{dt} =\frac{2\mu}{\eps} \cosh(\overline{V}) \big(r_1( t)-r_2(t)-W)\big)
\feq
\beq \label{W2}
\frac{dW}{dt} =\frac{2\mu}{\eps} \sqrt{1+ \sinh^2(\overline{V})} \big(r_1( t)-r_2(t)-W\big)
\feq

\beq \label{W3}
\frac{dW}{dt} =\frac{2\mu}{\eps} \sqrt{  1+( \sinh(\sinh^{-1}(W/\mu))^2     } \big(r_1( t)-r_2(t)-W\big)
\feq

\beq \label{W4}
\frac{dW}{dt} =\frac{1}{\eps} \sqrt{  4\mu^2+W^2    } \big(r_1( t)-r_2(t)-W\big)
\feq
et : 
$$\Delta(\eps,\mu) = \frac{1}{2\pi} \int_0^{2\pi}\displaystyle  r_1( s)+r_2(s) + 2\mu\left(  \sqrt{1+(W_{\mu}/2\mu)^2}-1  \right) ds $$
où $W_{\mu}$ est la solution périodique de \eqref{W4}. Si on pose $m = 2\mu$ et $W = y$ il vient finalement :
\beq \label{W5}
\frac{dy}{dt} =\frac{1}{\eps} \sqrt{  m^2+y^2    } \big(r_1( t)-r_2(t)-y\big)
\feq
qui, si l'on rajoute $\frac{dx}{dt} = 1$, est le système $S_m$ avec $f = r_1-r_2$. On a donc :
$$\Delta(\eps,\mu ) = \Delta(\eps,m) = \frac{1}{2\pi} \int_0^{2\pi}  r_1( s)+r_2(s) + \left(  \sqrt{m^2+y_m^2(s)}-m  \right) ds $$
où $y_m$ est la solution périodique de \eqref{W5}.\\\\
0n voit que si $m =  \sim 0$ on a :
$$\Delta(\eps,m) \sim  \frac{1}{2\pi} \int_0^{2\pi}\displaystyle  r_1( s)+r_2(s) + |y_m(s)| ds $$
La condition $\chi > 0$ de \eqref{HK} peut se lire :
$$2 \chi =  \frac{1}{2\pi} \int_0^{2\pi} 2 \max_{i = 1,2} r_i( s)ds = \frac{1}{2\pi} \int_0^{2\pi}r_1(s)+r_2(s) + |r_1(s)-r_2(s)| ds>0 $$
Donc démontrer la proposition \ref{conjecturekat} revient à montrer que, si l'on pose $m = \emat^{\rho/\eps} $  la quantité :
$$ \frac{1}{2\pi} \int_0^{2\pi}  |r_1(s)-r_2(s)| - y_m(s)  ds$$
peut être rendue  plus petite que la quantité standard $\chi > 0$ pour des valeurs non infiniment petites de $\rho$.

\paragraph{Analyse de la solution périodique $y_m$.} Pour fixer les idées nous le faisons pour un système particulier ; l'analyse générale procède de façon évidente de la même manière.

On considère le système :
\beq \label{per1}
S_m \quad \quad \left \{
\begin{array}{lcl}
\displaystyle \frac{dx}{dt} &=&1\\[8pt]
\displaystyle\frac{dy}{dt}& =&\displaystyle \frac{1}{\eps} \sqrt{  m^2+y^2    } \big(r_1( t)-r_2(t)-y\big)
\end{array}
\right.
\feq
avec $r_1-r_2$ de période $2\pi$ défini par : :
\beq \label{r1mr2}
\left[  
\begin{array}{lcl}
\displaystyle x \in [0,\pi[ &\Longrightarrow &r_1(x)-r_2(x) = \cos(x) \\[4pt]
\displaystyle x \in [\pi,2 \pi[& \Longrightarrow &r_1(x)-r_2(x) = 1 + 1.5 \sin(x)
\end{array}
\right.
\feq
Cette fonction est continue en tout point, sauf $x = \pi$, elle change de signe pour les quatre valeurs $\theta_1, \theta_2, \theta_3, \theta_4$ (voir figure \ref{grapher1r2}).

Pour $m >0$,  $S_m$ possède une unique solution périodique globalement asymptotiquement stable $y_m$ que nous analysons maintenant pour des valeurs exponentiellement petites de $m$ de la forme $m = \emat^{\rho/\eps}\,\;\rho \lnsim 0$ .  L'espace des phases de \eqref{per1} est $S^1\times \Rmat$ où $S^1 = \Rmat (\mathrm{mod} 2\pi)$. 
  \begin{figure}[ht]
  \begin{center}
 \includegraphics[width=1.0\textwidth]{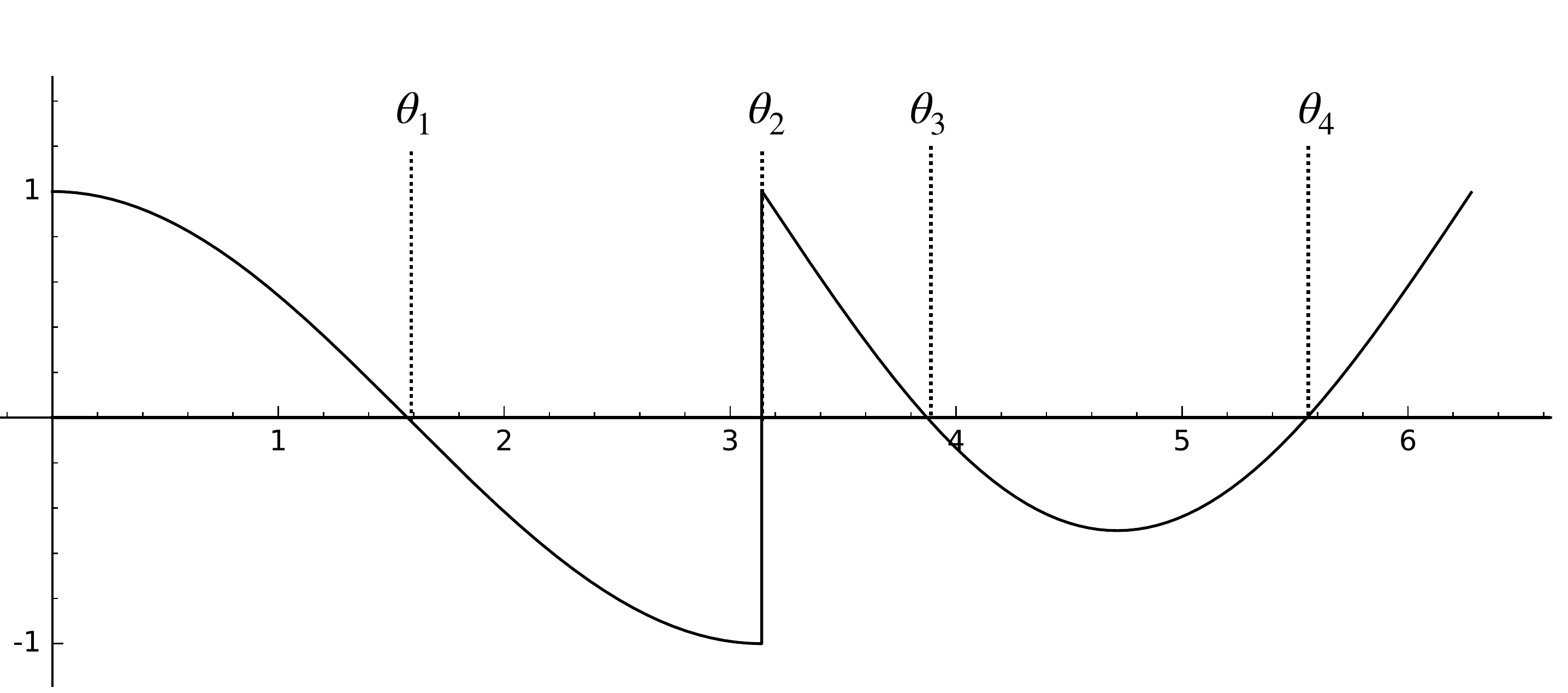}
 \caption{Graphe de $r_1-r_2$} \label{grapher1r2}
 \end{center}
 \end{figure}

Au système $S_m$ on associe le système contraint :
\beq \label{Cper1}
S_m^0 \quad \quad \left \{
\begin{array}{lcl}
\displaystyle \frac{dx}{dt} &=&1\\[8pt]
\displaystyle 0 & =&\displaystyle \sqrt{  m^2+y^2    } \big(r_1( t)-r_2(t)-y\big)
\end{array}
\right.
\feq

\begin{lemme}
Il existe $\rho^* \lnsim 0$ tel que pour tout $\rho \geq \rho^*$ on a : $$S(\theta_i)  < \theta_{i+1}\quad \quad (\mathrm{mod}\,2\pi)$$
\end{lemme}
\textbf{Preuve} C'est une conséquence immédiate de la définition de $S(x)$ (voir proposition \ref{rderho})\\\\
\noindent A partir de maintenant nous supposons que :
$$ \rho^*\leq \rho \lnsim 0$$
  \begin{figure}[ht]
  \begin{center}
 \includegraphics[width=1.0\textwidth]{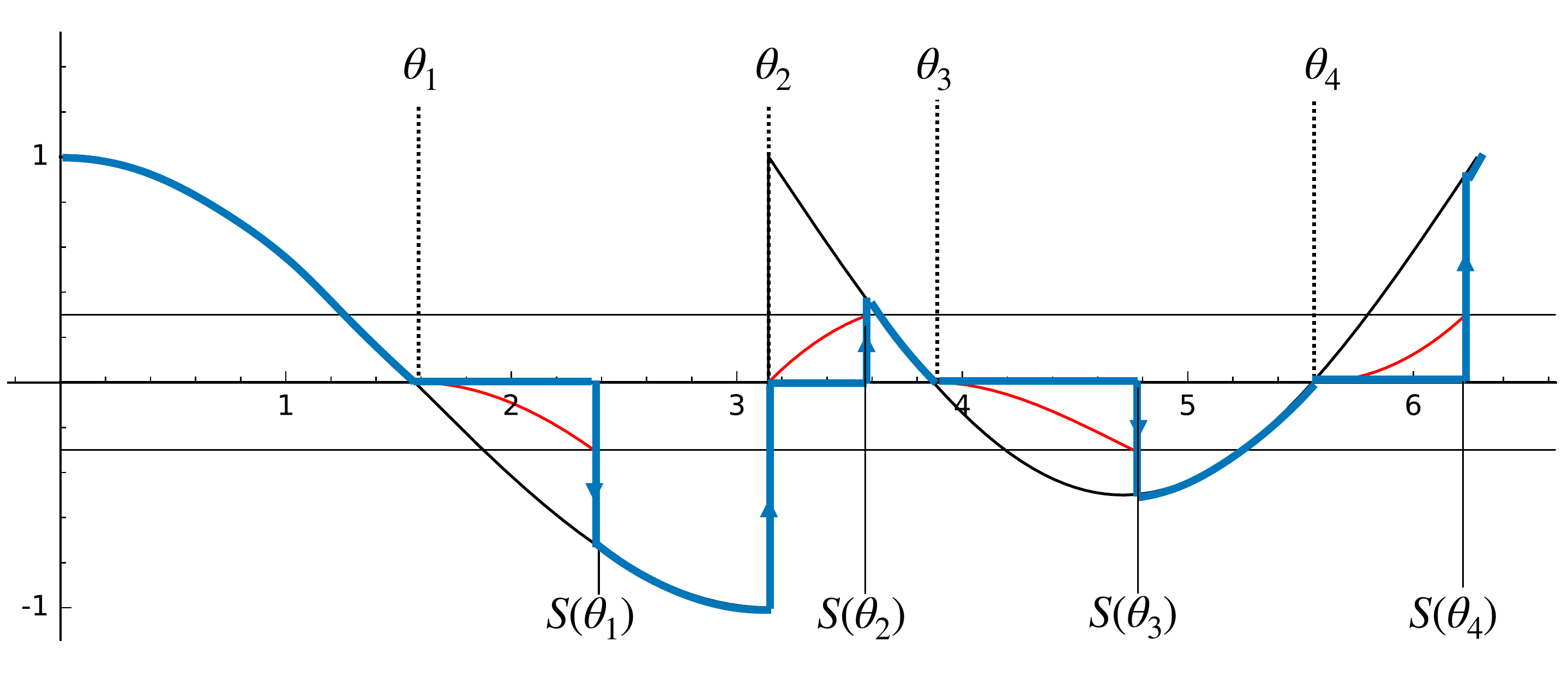}
 \caption{C-trajectoire périodique.} \label{grapher1r2bis}
 \end{center}
 \end{figure}
   \begin{figure}[ht]
  \begin{center}
 \includegraphics[width=1.0\textwidth]{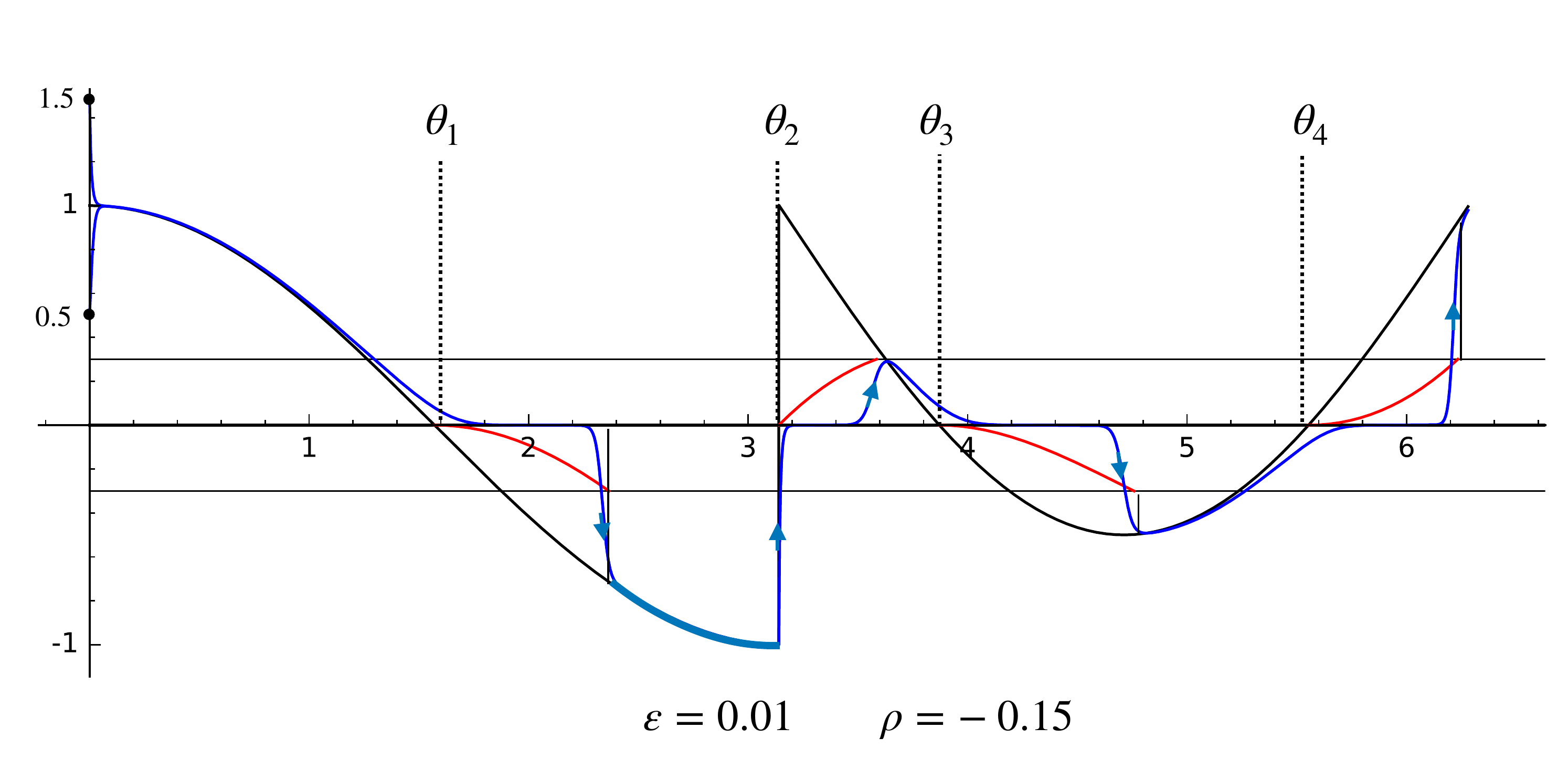}
 \caption{deux trajectoires de $S_m$ qui encadrent la C-trajectoire périodique.} \label{grapher1r2bister}
 \end{center}
 \end{figure}
La figure \ref{grapher1r2bis} représente la C-trajectoire issue du point : $(0,1)$. Elle est périodique. Sur la figure \ref{grapher1r2bister} nous avons les solutions issues respectivement des conditions initiale $(0,0.5)$ et $(0,1.5)$. Nous pouvons leur appliquer le théorème d'approximation et donc, après avoir longé respectivement les segments $[(0,5), (0,1)]$ et $[(0,1.5), (0,1)]$ ces deux trajectoires longent la C-trajectoire périodique ; ceci montre que l'image du segment $[(0.5),(0,1.5)]$ par l'application de Poincaré est contenue dans le halo de $(0,1)$ donc que la solution périodique de $S_m$ est infiniment proche de la C-trajectoire périodique. Puisque $\max{S(\theta_i) - \theta_i}$ tend vers $0$ quand $|\rho|$ tend vers $0$ ceci achève la preuve de la proposition \ref{conjecturekat}. En effet, c'est ici qu'intervient de façon essentiel le fait que  $r_1,r_2$ sont des fonctions standard ce qui entraine que $\chi >0$ est standard, donc  que $\chi \gnsim 0$ et, comme les $S(\theta_i)$ sont continues et standard, il existe un $\eta < 0$ standard tel que :
$$\rho \leq \eta \Longrightarrow \int_\Omega (r_1(s)-r_2(s)) ds <\frac{ \chi}{2}$$
avec $\Omega = \cup_i[\theta_i,S(\theta_i)]$

\newpage
{\huge\textbf{Appendice}}
 \appendix
 \section{Introduction à l'ANS}\label{ANS}
 
On attribue l'invention de l'ANS. à A. Robinson \cite{ROB16}. Il en existe un grand nombre de versions et on peut consulter l'article de revue \cite{FLE17} pour en avoir une idée. La version sur laquelle je m'appuie ici, appelée {\em Internal Set Theory},  (I.S.T.) est due à E. Nelson. Il l'a proposée dans l'article \cite{NEL77} qui contient également une introduction à la pratique de I.S.T. 
Plus de quarante ans après la publication de cet article il faut bien constater que l'espoir de Nelson et de ses adeptes,  de voir la pratique de l'A.N.S. se généraliser a été déçu. Peut être que le prix du ''ticket d'entrée'' était  trop élevé.

Par la suite Nelson a également proposé d'utiliser des versions plus élémentaires, certes moins puissantes, mais dont il a cependant montré toute la fécondité dans son livre {\em Radically Elementary Probability Theory} \cite{NEL87}. A l'instar  de ce que que  son ami Reeb proposait (\cite{REE79} (voir aussi  \cite{LUT87,LUT92}) il a constamment encouragé les mathématiciens à utiliser des théories faibles de l'ANS comme par exemple dans cet extrait de son article {\em The virtue of simplicity} \cite{NEL07}. \dcom
Much of mathematics is intrinsically complex, and there is delight to be found in mastering complexity. But there can also be an extrinsic complexity arising from unnecessarily complicated ways of expressing intuitive mathematical ideas. Heretofore nonstandard analysis has been used primarily to simplify proofs of theorems. But it can also be used to simplify theories. There are several reasons for doing this. First and foremost is the aesthetic impulse, to create beauty. Second and very important is our obligation to the larger scientific community, to make our theories more accessible to those who need to use them. To simplify theories we need to {\em  have the courage to leave results in simple, external form}\footnote{Souligné par moi}—fully to embrace nonstandard analysis as a new paradigm for mathematics.
\fcom
J'ai essayé de me conformer à cette recommandation que je discute un peu plus loin dans l'appendice.

  \subsection{Etre standard ou non, interne ou non : le langage étendu.}\label{tutoANS}
 Pour qui ne possède pas des notions élémentaires de logique, ce qui est le cas de beaucoup de mathématiciens, la pratique de la théorie I.S.T. de Nelson est un peu déroutante au début. Rappelons quelques notions élémentaires de logique sans lesquelles on ne peut pas comprendre la démarche de I.S.T..
 
  La {\em théorie formelle des ensembles} Z.F.C. (voir le classique \cite{KRI69} pour un exposé très accessible) est la théorie à laquelle les textes mathématiques contemporains, dans leur quasi totalité,  sont supposés se référer. 
  Rappelons que dans cette théorie les formules comportent en plus des symboles logiques ($\forall,\;\exists,\;=,\cdots$ ) un prédicat binaire $x\in y$ qui se lit ''$x$ appartient à $y$''. Z.F.C. propose une série d'axiomes qui permettent de construire de proche en proche, à partir d'un ensemble primitif, l'ensemble vide, de nouveaux ensembles qui portent des noms familiers : $\Nmat, \mathbb{Z}, \mathbb{Q}, \Rmat, \cdots$. 
  
  Les fonctions aussi sont considérées comme étant des ensembles. Par  exemple :
 $$(x \mapsto x^2 )\stackrel{\mathrm{def}}{\longleftrightarrow} \{(x,y)\in \Rmat^2\;;\; y = x^2\}$$
 Dans l'expression ci-dessus, à  la suite du   "$\;;\;$",  à l'intérieur des crochets  $\{\}$ on a l'écriture "$y = x^2$". C'est une ''formule'' (un ''énoncé'') $\phi(x,y)$ qui est ou n'est pas satisfaite par le couple $(x,y)$. Le fait que $ \{(x,y)\in \Rmat^2 \;;\;  y = x^2\}$ découle de l'axiome :
 \bitbul
 \item  $A_c \quad \quad$ Si  $\phi$ est une "formule" qui ne comporte que "$z$" comme variable libre (i.e. non quantifiée) :$ \quad \forall x\; \exists y \; \forall z \;:\;  (z \in y  \Leftrightarrow z \in y \;\mathrm{et} \; \phi(z))$
 \fit
 Une ''formule'' n'est pas n'importe quoi. C'est une suite de parenthèses, de symboles logiques, d'occurrences   du symbole d'appartenance $\in$ et de constantes (les noms des ensembles déjà construits) qui respecte une syntaxe de construction ; les formules correctement écrites constituent le langage $\mathcal{L}$ de Z.F.C. 

\paragraph{ I.S.T. est une extension de Z.F.C.} On augmente les possibilités d'expression de $\mathcal{L}$ en ajoutant le prédicat $st(x)$ qui se lit "$x$ {\em est standard}" ; les règles de manipulations de $st$ sont définies par trois axiomes, que je ne précise pas tout de suite\footnote{D'ailleurs, à part les spécialistes, quels mathématiciens ont présent à l'esprit, lorsqu'ils rédigent leurs démonstrations,  la dizaine d'axiomes de Z.F.C. ?}, d'où se déduisent rapidement des propriétés utiles de $st$. Par exemple on a :
 \ben
 \item $st(0)\stackrel{\mathrm{se \;lit}}{ \longleftrightarrow}$ $0$ est standard.
 \item $n \in \Nmat : st(n) \Rightarrow st(n+1)\stackrel{\mathrm{se \;lit}}{ \longleftrightarrow}$ le successeur d'un standard est standard.
 \item $\exists \; \omega \in \Nmat \;\forall n \in \Nmat : (st(n) \Rightarrow n < \omega)  \stackrel{\mathrm{se \;lit}}{ \longleftrightarrow}$ il existe un entier infiniment grand (plus grand que tous les entiers standard).
 \fen
 On peut être tenté de voir ici un problème: un théorème classique de Z.F.C. dit que si $E$ est une partie de $\Nmat$, si $0\in E$ et si $n \in E \Rightarrow (n+1) \in E$ alors $E = \Nmat$. Les points  1. + 2. ci-dessus n'entrainent-ils pas que l'ensembles des entiers standard est $\Nmat$ tout entier? Ce qui contredit 3. En réalité il n'y a pas de problème car rien ne nous permet d'affirmer que les entiers standard constituent un ensemble :{\em  $\;st(x)$ n'étant pas une formule du langage $\mathcal{L}$ de Z.F.C. l'axiome $A_c$ ne peut pas être invoqué.} A partir de l'existence d'entiers infiniment grands on définit toutes les notions et symboles utiles de l'analyse élémentaire tels que {\em infiniment petit} ou $\sim$ : 
 $$x \sim y \Longleftrightarrow |x-y| \;\mathrm{infiniment\;petit}$$
qui sont des notions {\em externes}, c'est à dire externes au langage traditionnel de Z.F.C.
 
 C'est pourquoi partir, à partir  de maintenant, nous devons avoir présent à l'esprit  qu'il y a deux sortes de formules\;;\;
 \bito 
 \item  Les {\em formules internes} $\phi$, celles qui appartiennent au langage $\mathcal{L}$ de Z.F.C. pour lesquelles on peut écrire grace à l'axiome de compréhension :
 $$ E = \{x \in y \;;\; \phi(x)\}$$ 
 où $E$ est un ensemble, 
 \item  les formules $\Phi$ qui ne sont pas internes pour lesquelles l'écriture :
 $$ \{x \in y \;;\; \Phi(x)\}$$ 
 ne désigne pas un ensemble ; pour éviter le risque de croire qu'on affaire à un ensemble, quand $\Phi$ n'est pas une formule interne j'écrirai : 
 $$ \eset{x \in y \;;\; \Phi(x)}$$
 \fit
 Par exemple, dans $\Rmat^2$ si $\mathcal{C}$ est le graphe d'une fonction $f$  nous pourrions écrire :
 $$ \mathrm{Halo}(\mathcal{C}) = \eset{ (x,y) \in \Rmat^2   \,;\, y \sim f(x)}$$ 
 Ici la formule $\Phi$, qui porte sur les couples $(x,y)$, est la formule externe $y \sim f(x)$.
 
 Il faut bien faire attention à distinguer :
 \bitbul
 \item Etre standard  qui  est une qualité que peut avoir ou pas  un ensemble (un objet mathématique).
 \item Etre interne ou non qui  est une qualité que peut avoir ou pas  une formule.
 \fit
 
 On démontre que I.S.T. est une {\em extension conservative de Z.F.C.} ce qui veut dire que si A et B sont deux propositions internes (qui s'expriment dans $\mathcal{L}$) et si il existe une preuve dans I.S.T. de $A \Rightarrow B$  alors il en existe une preuve dans $\mathcal{L}$. 
 L'importance de la conservativité réside dans le fait que si I.S.T. était contradictoire alors Z.F.C. le serait aussi. En effet, si I.S.T. était contradictoire il y existerait une preuve (par exemple) de la proposition interne $1 =2$, et, comme I.S.T. est conservatif, il existerait une preuve dans Z.F.C. de $1 = 2$, donc une contradiction dans Z.F.C. Du point de vue du risque de contradiction, donc de la rigueur, il n'est pas plus dangereux de travailler dans I.S.T. que dans Z.F.C.

  \subsection{Qui est standard}
  Les axiomes de I.S.T. sont les règles de manipulation du prédicat $st(x)$. Les voici. 
  On utilise les abréviations ~:\\
$\forall^{st}x \longleftrightarrow \forall x \;(x\; \mathrm{standard}) \Rightarrow$\\
$\forall^{fin}x \longleftrightarrow \forall x \;(x\; \mathrm{fini}) \Rightarrow \quad$\\
$\forall ^{st\;fin}\longleftrightarrow \forall^{st} x \;(x\; \mathrm{fini}) \Rightarrow $\\
$ \exists^{st} x \longleftrightarrow \exists x\; (x\; \mathrm{standard})\; \wedge\;$\\
Le mot {\em fini} est pris dans son sens mathématique usuel~: un ensemble $x$ est {\em fini } ssi toute injection de $x$ sur lui même est surjective.
\bitbul
\item {\em  \textbf{I}déalisation}~: Soit $B(x,y)$ une formule {\em interne} :
$$\forall ^{st\;fin} z\, \exists\, x\forall \,y \in z \; B(x,y) \Longleftrightarrow  \exists \,x \;\forall^{st}y\; B(x,y)$$
\item {\em \textbf{S}tandardisation}~: Soit $\phi(z)$ une formule, interne ou non :: 
$$\forall ^{st}x\,\exists^{st}y\,\forall^{st}z\;(z \in y \Leftrightarrow z\in x \wedge \phi(z))$$
\item {\em \textbf{T}ransfert }: Soit $A(x,t_1,\cdots,t_k)$ une formule {\em interne} sans autre variable libre que $(x,t_1,\cdots,t_k)$. Alors~:
$$\forall ^{st}t_1 \cdots \forall ^{st}t_k \; \big( \forall ^{st}x\;A(x,t_1,\cdots,t_k) \Rightarrow \forall x \;A(x,t_1,\cdots,t_k) \big)$$
et par contraposition~:
$$\forall ^{st}t_1 \cdots \forall ^{st}t_k \; \big( \exists \; x\;A(x,t_1,\cdots,t_k) \Rightarrow \exists ^{st} x \;A(x,t_1,\cdots,t_k) \big)$$
\fit
d'où le nom donné à la théorie : \textbf{I.S.T.}, acronyme qui vaut aussi pour {\em Internal Set Theory}.

On remarque la similitude entre l'axiome $A_c$ de Z.F.C. et l'axiome des standardisation. Le premier nous dit que pour toute formule, les éléments d'un ensemble $x$ qui satisfont cette formule constituent un ensemble, le second précise que si l'ensemble de départ est standard, l'ensemble dont les éléments standard satisfont la formule peut être choisi standard. Comme pour l'essentiel les nouveaux ensembles de Z.F.C. sont construits par l'application de $A_c$ il s'ensuit que tous les objets mathématiques (les ''ensembles'') construits de façon unique sont standard.  Ainsi tous les objets des mathématiques traditionnelles tels que $\Nmat, \Rmat, x\mapsto  \cos (x),  \cdots$ sont standard. Mais tous les éléments d'un ensemble standard ne sont pas forcément standard. En fait tous les éléments d'un ensemble standard sont standard si et seulement si ce dernier est un ensemble fini. 

Quelques exemples éclairent la situation.
\bit
\item $\Rmat$ est standard, il contient des éléments non standard, les infinitésimaux, les infiniment grands, mais pas seulement. Par exemple $1+\eps$ avec $\eps$ i.p. est un non standard, ni infiniment grand, ni infiniment petit.
\item $\Rmat^2$, $\Rmat^3$ sont standard, $\Rmat^n$ est standard ssi $n$ est standard.
\fit
 La fonction de $\Rmat^2$ dans $\Rmat$ :
$$(x, m) \mapsto \frac{\emat^{mx}-\emat^{-mx} }{\emat^{mx}+\emat^{-mx}} = \tanh(mx)$$ 
 est une ''vraie'' fonction au sens traditionnel ; quels que soient $m$ et $x$ qu'ils soient  standard ou non le réel $\emat^{mx}$ est bien défini de façon interne par des formules (ici la somme d'une série convergente) qui n'utilisent pas le prédicat $st()$. C'est une fonction standard du couple $m, x$. Mais si maintenant $m$ est compris comme un paramètre, pour tout $m$ dans $\Rmat$  :
 $$x \mapsto f_m(x) = \mathrm{th}(mx)$$
 est une fonction de $\Rmat$ dans $\Rmat$ qui, en tant que fonction, est standard si  $m$ est standard et non standard si $m$ ne l'est pas. Prenons $m$ i.g., donc $\mathrm{th}(mx)$ est non standard, mais elle n'a pas cessé d'être la fonction traditionnelle $\mathrm{th}(mx)$ qui est continue au point $x = 0$. Toutefois dès que $x$ est non i.p., le produit $mx$ est i.g. (négatif ou positif) et la valeur prise par $\mathrm{th}(mx)$ est infiniment proche de $±1$. Lorsque $m$ est i.g. la fonction $f_m$ a la particularité qu'un accroissement infinitésimal de la variable $x$ au voisinage de $0$ peut entrainer un accroissement non infinitésimal de $f$. On dit que la fonction est "continue" (au sens traditionnel) mais n'est pas {\em S-continue} (on dit  qu' une fonction $f$ est {\em  S-continue} au point $x$ si $dx\sim 0$  entraine $f(x+dx) \sim f(x)$).
 
  Tout ceci n'est qu'une façon de parler dans le nouveau langage étendu de la convergence (non uniforme) de la famille de fonctions $f_m$ vers la fonction $signe(x)$ lorsque $m$ tend vers l'infini. Pour cet exemple le bénéfice
 ne semble pas évident mais il le devient dans le cas qui nous intéresse de perturbations singulières d'équations différentielles.
 
 En effet, on s'intéresse à une famille d'objets ''équation différentielle'' :
 \beq
\Sigma_{\eps} \quad \quad  \displaystyle \eps \frac{dx}{dt} = f(x) \Longleftrightarrow \frac{dx}{dt} = \frac{1}{\eps} f(x)
 \feq 
 pour $\eps$ tendant vers $0$. On aimerait démontrer des théorèmes du genre ''le portait de phase de $\Sigma_m$ tend vers celui de $\Sigma_0$'' mais on ne peut pas parce que $$0\, \frac{dx}{dt} = f(x)$$ n'est pas une équation différentielle. En ANS on fixe $\eps >0$ infiniment petit et on s'intéresse au portrait de phase de $\Sigma_{\eps}$ qui reste parfaitement défini et on cherchera à établir des résultats indépendants de la valeur de $\eps$, pourvu que cette dernière soit infiniment petite comme il l'est dans le théorème d'approximation des trajectoires de cet article.
 
 Un dernier mot sur les fonctions standard qui nous intéresse ici. Soit $x \mapsto f(x)$ une fonction standard à valeurs réelles dont la borne inférieure $\inf(f)$ est strictement positive ; alors $\inf(f) \gnsim 0$. En effet, la borne inférieure étant définie de façon unique à partir d'un objet standard est, à son tour, standard ; un nombre standard strictement positif n'est pas infinitésimal.
 
\subsection{Théories plus faibles que I.S.T.}\label{cheapANS}

Il est indéniable que les formules qui définissent les axiomes I.S. et T.  ont un effet repoussoir sur le mathématicien ordinaire. C'est pourquoi des mathématiciens comme G. Reeb \cite{REE79} et R. Lutz \cite{LUT87, LUT92} ont proposé de simplifier cette axiomatique en l'affaiblissant. Dans \cite{LUT87} il est proposé un ''sous système'' de I.S.T. appelé Z.F.L. (''L.'' en l'honneur de Leibnitz) dont Lutz nous dit dans \cite{LUT92} :
\dcom
Dans I.S.T. il existe des équivalences automatiques entre  théorèmes internes et formulations externes (voir Nelson \cite{NEL77}). Il n'en est pas de même dansZ.F.L ; des formulations concurrentes des idées intuitives de limite, continuité etc. cohabitent en restant irréductibles l'une à l'autre. 
\fcom
Nelson lui même a proposé une version affaiblie de I.S.T. dans son livre {\em Radically Elementary Probability Theory} où il introduit un appendice montrant l'équivalence formelle des résultas radicalement élémentaires avec les résultats traditionnels par les remarques :
\dcom
The purpose of this appendix is to demonstrate that theorems of the conventional theory of stochastic processes can be derived from their elementary analogues by arguments of the type usually described as generalized nonsenses; there is no probabilistic reasoning in this appendix. This shows that the elementary nonstandard theory of stochastic processes can be used to derive conventional results; on the other and, it shows that neither the elaborate machinery of the conventional theory nor the devices from the full theory of nonstandard analysis, needed to prove the equivalence of the elementary results with their conventional forms, add anything of significance : the elementary theory has the same scientific content as the the conventional theory. This is intended as a self-destructing appendix.
\fcom

Ces auteurs nous proposent un langage externe simplifié où ne subsiste que la structure des nombres réels enrichie par la possibilité de parler d'ordre de grandeur et nous encouragent à y travailler ; le prix de cette simplicité est l'impossibilité de traduire dans le langage conventionnel ce qui est dit de façon externe. Mais est-ce bien un prix à payer ? N'est-il pas, au contraire, plus enrichissant de garder les deux points de vues ? Nous avons remarqué que pour $m$ infiniment grand la fonction $m \mapsto \tanh (mx)$ est (traditionnellement) continue, mais pas {\em S-continue}. En revanche la fonction $x \mapsto \eps [x/ \eps]$ (où $[\cdot]$ désigne la partie entière), lorsque $\eps \sim 0$, est S-continue, mais n'est pas continue : elle progresse par petits sauts infinitésimaux et ne satisfait pas au théorème de la valeur intermédiaire. La continuité et la S-continuité ne parlent pas exactement  de la même chose lorsqu'on les considère sur des objets non standard. 

D'autre part cette question du retour d'énoncés externes vers des énoncés internes à $Z.F.C.$ est liée à l'utilisation de l'axiome du choix. Une critique souvent formulée contre l'usage de l'ANS est son caractère essentiellement ''non constructif''. Il est vrai que la démonstration de l'existence d'un modèle de l'ANS ayant une force au moins égale à I.S.T. à l'intérieur de $Z.F.C.$ repose sur l'axiome du choix dans sa version la plus générale. Mais ce n'est pas vrai pour les versions faibles. Les logiciens ont montré que, à condition de renoncer à transformer les énoncés externes en énoncés internes des versions axiomatiques faibles convenables permettent 
un pratique nonstantard de l'essentiel des mathématiques en n'utilisant, au plus, que l'axiome du choix dénombrable (voir \cite{HRB21} pour un exposé relativement accessible).

C'est pour toutes ces raisons que j'ai essayé de ne pas utiliser toute la puissance de I.S.T., en particulier l'axiome $S.$ qui permet de définir {\em l'ombre} d'une fonction comme le {\em standardisé} du {\em halo}  de son graphe, ce qui aurait obligé le lecteur à comprendre ce qu'est le standardisé, et donc aurait demandé un effort supplémentaire. 
En revanche je n'ai pas cherché à définir précisément le système axiomatique nonstandard  permettant de formaliser au mieux l'intuition naturelle des praticiens de la théorie des équations différentielles singulièrement perturbées. C'est un travail qui reste à faire.

\section{ANS et systèmes dynamiques} \label{ANSetEDO}

 \subsection{Dépendance continue des solutions} \label{gronwall}
 Le  très classique théorème de dépendance continue des solutions d'une E.D.O. par rapport aux conditions initiales et aux perturbations que l'on déduit du  lemme de Gronwall possède l'expression externe suivante :
 \begin{theoreme} Soit $ f $ $C^1$-limitée et globalement lipschitz, de constante de lipschitz $k$ limitée. Soit $g(x)$ de classe $C_1$
 infiniment petite (i.e. pour tout $x$ on a $g(x) \sim 0$). Soient $x_0\sim y_0$  et soient $x(t)$ et $y(t)$ les solutions respectives de :\\
 $\frac{dx}{dt} = f(x);\quad x(0) = x_0$\\
$\frac{dy}{dt} = f(y)+g(y);\quad y(0) = y_0$\\
 Alors, pour tout $T$ non infiniment grand :
 $$t \in [0,T]   \Longrightarrow y(t) \sim x(t)$$
 \end{theoreme}
 \textbf{Preuve} Soit $E$ un ensemble compact contenant la trajectoire $\{ x(t) ; t \in [0,T]\}$ dans son intérieur. Le maximum de $|g(x)|$ sur $E$ est $m \sim 0$. L'inégalité de Gronwall classique dit que, tant que $y(t) \in E $ on a :
 $$ |x(t)-y(t)| \leq |x_0-y_0|\emat^{kt m} \sim 0$$
comme $k$, $t$ et $m$ sont limité il en est de même de $\emat^{kt m}$. Donc $y(t)$  ne peut pas sortie de $E$ avant $T$ et $y(t) \sim x(t)$ jusqu'à $T$.\\
 $\Box$\\
 Nous laissons au lecteur le soin de formuler le cas vectoriel et de s'affranchir de l'hypothèse de l'existence d'une constante de Lipschitz globale pour $f$.  
 \subsection{Entrée dans le {halo} d'une courbe lente.}\label{halocourbelente}
 Commençons par rappeler un résultat élémentaire bien connu concernant les champs de vecteurs du plan. On considère un contour $\Gamma$, fermé, sans point double,  \\
 \begin{minipage}{0.55 \textwidth} 
 différentiable par morceaux, constitué de la réunion de segments $\arc{ab}$,  $\arc{bc}\cdots \arc{ga}$ (figure ci-contre), qui définit un domaine $\Gamma$ et un champ de vecteur différentiable $X$, {\em qui ne s'annule pas}  dans $\Gamma \cup \partial \Gamma$ et qui pointe strictement vers l'intérieur de $\Gamma$, sauf pour un des segments ($\arc{ef}$ sur la figure) où il pointe vers l'extérieur. 
 \end{minipage}
 \begin{minipage}{0.40 \textwidth}
 \includegraphics[width=1\textwidth]{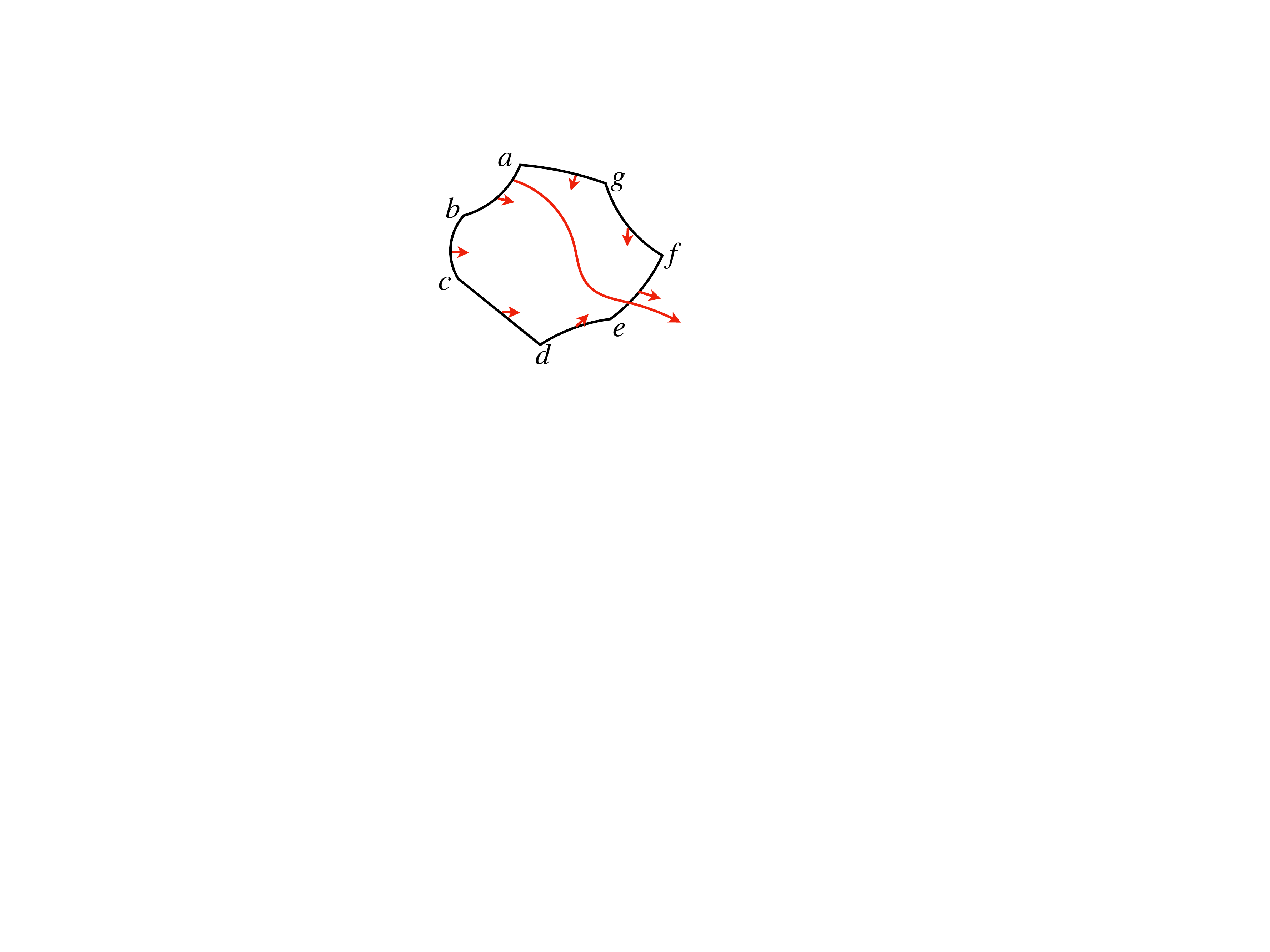}
 \end{minipage}
 \begin{theoreme}\label{piege}
 Pour toute condition initiale dans $\Gamma \cup \partial \Gamma$ la trajectoire de $X$ qui en est issue sort de $\Gamma$ en un point de $\arc{ef}$.
 \end{theoreme} 
  On considère maintenant  le système :
  \beq \label{eqLR3}
  \left\{
\begin{array}{lcl}
\displaystyle \frac{dx}{dt}& =& 1\\[6pt]
\displaystyle  \frac{dy}{dt} &=& \displaystyle \frac{1}{\eps}\Big(f(x)-y\Big)
 \end{array} 
 \right.
\feq
avec $f$ $C^1$-limitée. On suppose $\eps \sim 0$  fixé. 0n dit qu'un tel système est {\em lent-rapide}. Le graphe de $f$ s'appelle la {\em courbe lente}.

Soit $(x_0,y_0)$ une condition initiale limitée à l'instant $0$, n'appartenant pas au halo $\mathcal{H}$ de la courbe lente (au dessus du graphe de $f$ pour fixer les idées).
\begin{proposition}\label{cl1}
Soit $(x(t),y(t))$ la solution de \eqref{eqLR3} de condition initiale $(x_0,y_0)$. Il existe $t^*>0$   :
\ben
\item $t^*\sim 0$
\item $t \in \,[0,t^*] \Longrightarrow x(t) \sim x_0$
\item $ y(t^*) \sim f(x_0)$
\item $ t > t^*$ et {\em  limité} $ \Longrightarrow y(t) \sim f(x(t))$
\fen
En d'autres termes plus géométriques : la trajectoire longe en descendant le segment $\{(x_0,y) : y \in [f(x_0),y_0]\}$, pénètre dans le halo du graphe de $f$ avant le temps $t^* \sim 0$   et y reste tant que $t$ est limité (au delà on ne peut rien dire comme le montre l'exemple $f(x) = x^2$ qui s'intègre explicitement). 
\end{proposition}
\textbf{Preuve.}\\
On suppose, pour fixer les idées,  que $f'(x_0) < 0$.\\
 Soit $\alpha \gnsim 0$ telle que le graphe de $f(x)+\alpha$ soit en dessous de $y_0$ ; il en existe puisque $(x_0,y_0)$ n'est pas dans le halo du graphe de $f$. Soit $\theta \gnsim 0$. On considère le contour $a,b,c,d$ défini par :\\
  \begin{minipage}{0.50 \textwidth} 
  \bito 
 \item Segment $\arc{a,b}$ : segment de droite joignant $a$ à $b$ ; $a$ est la condition initiale $(x_0,y_0)$, $b$ est le point $(x_0,f(x_0))$ d'intersection de la verticale $x = x_0$  avec le graphe de $f$.
  \item Segment $\arc{b,c}$ : portion du graphe de $f$ comprise entre $b$ et $c$ ; $c$ est le point d'intersection du graphe de $f$ avec la verticale $x = x_0+\theta$.
  \fit
\end{minipage}
  \begin{minipage}{0.5 \textwidth}
 \includegraphics[width=1\textwidth]{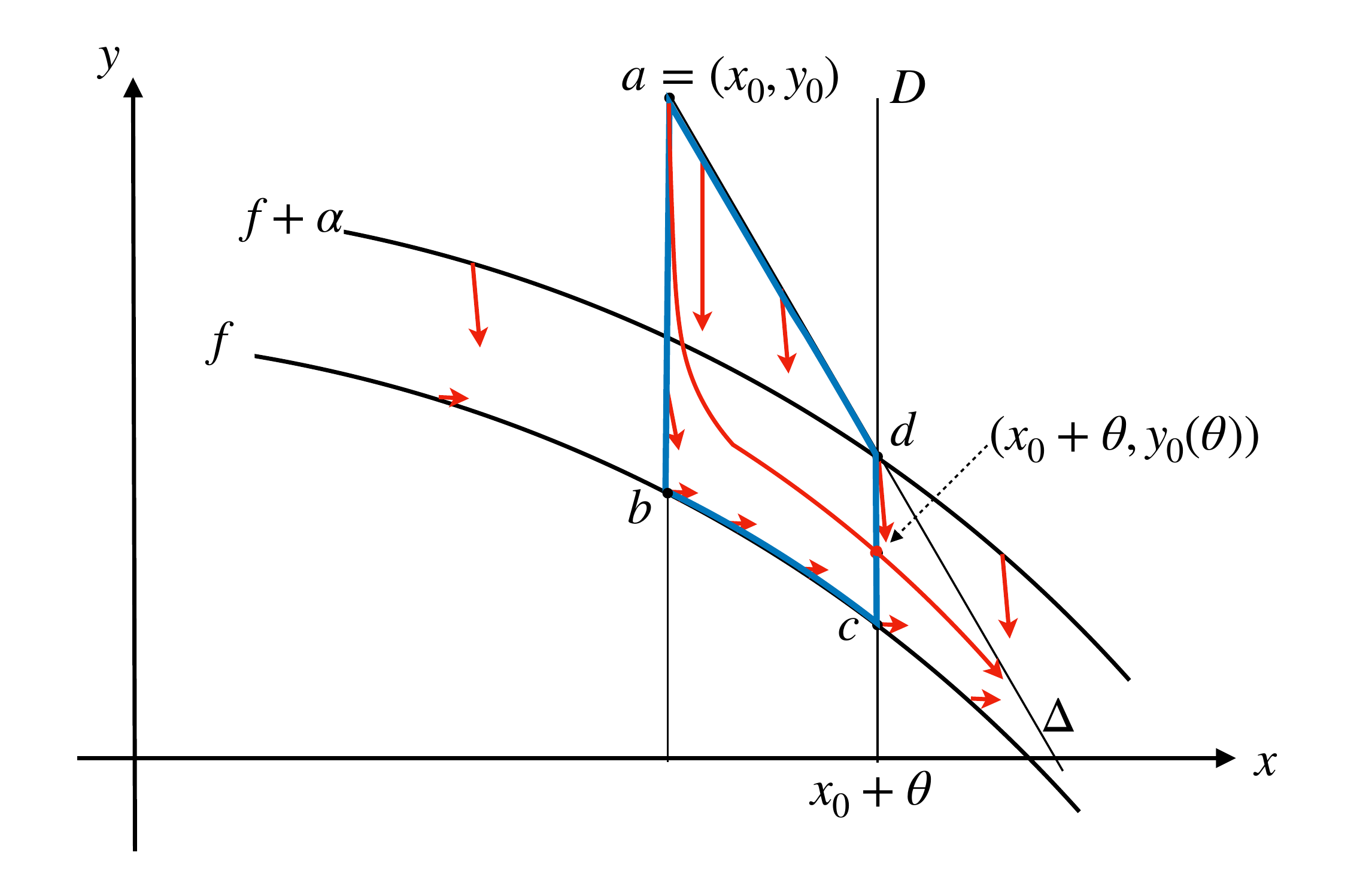}
 \end{minipage}
 \bito 
  \item Segment $\arc{c,d}$ : segment de droite joignant $c$ à $d$ ;  $d$ est le point d'intersection du graphe de $f+\alpha $ avec la verticale $x = x_0+\theta$.
 \item Segment $\arc{da}$ : segment de droite joignant $d$ à $a$.
 \fit 
Sur tous ces segments, sauf le segment $\arc{c,d}$, le champ de vecteur défini par \eqref{eqLR3}) pointe strictement vers l'intérieur de $\Gamma$ ; pour le segment $\arc{d,a}$ cela est une conséquence de ce que pour $(x,y)$ au dessus du graphe de $f(x)+\alpha$ on a $\frac{dy}{dx} <  \frac{ f(x)-y}{\eps}$ qui est {i.g.} (négatif), donc plus petit que la pente de $\Delta$. qui, elle, n'est pas infiniment grande parce que $\theta \nsim 0$.
Ceci montre que pour tout $\alpha \gnsim 0$ on a $f(x_0+\theta ) < y(\theta) < f(x_0+\theta)+\alpha)$ donc que $y(\theta)$ est infiniment proche de $f(x_0+\theta)$, donc que $(x(\theta),y(\theta))$ est dans le  {halo} du graphe de $f$. Nous laissons au lecteur le soin d'effectuer les modification évidentes qui permettent d'obtenir le même résultat dans le cas où la dérivée de $f$ en $x_0$ est positive ou nulle.

Nous avons démontré que pour tout $t\gnsim 0$ (aussi petit soit-il, mais non infiniment petit)  la solution $(x(t), y(t)) $ est dans le {halo} du graphe de $f$, mais la proposition dit  plus : qu'il existe un $t^*\sim 0$ tel que $\forall t \,\geq t^* \;:\; y(t)\sim f(y(t))$. Ce dernier point s'obtient à l'aide d'un argument typiquement ''non standard'' (appelé lemme de Robinson) qui est le suivant.
On considère l'ensemble :
$$
E =\left \{ t \in [0,1] \;:\; \frac{ |y(t)-f(x(t))|}{t} \leq 1 \right \}
$$
La définition de $E$ ne fait intervenir aucun terme du nouveau langage, c'est donc un ''véritable'' ensemble de la théorie des ensembles usuelle. Cet ensemble contient tous les réels non {i.p.} de $[0;1]$ ; en effet, pour $t\nsim 0$ nous avons $|y(t)-f(x(t)| \sim 0$ ce qui entraine que $\frac{ |y(t)-f(x(t))|}{t} \leq 1$ (on divise un  {i.p.} par un non  {i.p.}, le résultat est  {i.p.}).    La borne inférieure $l $ de $E$ est infiniment petite car si ce n'était pas le cas, $l/2$ ne serait pas infiniment petit et donc $l/2$                          
serait dans $E$ ce qui contredit la définition de la borne inférieure. Si la borne inférieure de $E$ est infiniment petite, $E$ contient nécessairement un {i.p.} $t^*$. Puisque $t^*$ est dans $E$ on a :
$$ \frac{|y(t^*)-f(x(t^*))|}{t^*} < 1$$
Pour que la quotient soit plus petit que 1 il faut que le numérateur soit {i.p.} puisque le dénominateur l'est.\\
$\Box$\\\\
Dans le résultat précédent il est essentiel que $|f'(x)|$ soit limité. Si ce n'était pas le cas le point $c$ pourrait être rejeté dans les infiniment grands négatifs et on ne pourrait plus garantir que le segment $\arc{da}$ est de pente limitée, donc que, le long  de $\arc{da}$, le champ pointe vers l'intérieur du domaine. 

Avec des adaptations évidentes de cette preuve on obtient  le résultat suivant un peu plus général que \eqref{eqLR3} qui peut être considéré comme la version ANS du théorème de Tichonov (en dimension 2) (voir \cite{TIK52, LOB98} :\\
Soit le système :
 \beq \label{eqLR4}
  \left\{
\begin{array}{lcl}
\displaystyle \frac{dx}{dt}& =& h(x,y)\\[6pt]
\displaystyle  \frac{dy}{dt} &=& \displaystyle \frac{1}{\eps}g(x,y)\Big(f(x)-y\Big)
 \end{array} 
 \right.
\feq
où $h$ est une fonction $C^1$-limitée qui ne s'annule pas.
\begin{proposition}\label{cl2}
Soit le système \eqref{eqLR4}. On suppose que :
$$\min g(x,y) \gnsim 0$$ 
Soit $(x(t),y(t))$ la solution de \eqref{eqLR3} de condition initiale $(x_0,y_0)$. Il existe $t^*>0 $  :
\ben
\item $t^* \sim 0$
\item $t \in \,[0,t^*] \Longrightarrow x(t) \sim x_0$
\item $ y(t^*) \sim f(x_0)$
\item $ t > t*$ et limité $ \Longrightarrow y(t) \sim f(x(t))$
\fen
\end{proposition}
\subsection{Systèmes semi-lent-rapides.}\label{demiLR}
Un champ semi-lent-rapide est un champ qui est infiniment grand d'un côté d'une courbe lente, non infiniment grand de l'autre côté. Je traite le cas particulier du champ semi-lent-rapide qui nous concerne. 

Soient les champs :
\beq \label{DLR1}
\Sigma \quad \quad  \left\{
\begin{array}{lcl}
\displaystyle \frac{dx}{dt}& =&1 \\[6pt]
\displaystyle \frac{dz}{dt} & =&\displaystyle z\sqrt{1+\varphi (z)  } ( f(x)-\psi(z))
\end{array} 
\right.
\feq
et
\beq \label{DLR2}
\tilde{\Sigma} \quad \quad  \left\{
\begin{array}{lcl}
\displaystyle \frac{dx}{dt}& =&1 \\[6pt]
\displaystyle \frac{dz}{dt} & =&\displaystyle zf(x)
\end{array} 
\right. \quad \quad  \quad \quad  \quad \quad  \quad \quad 
\feq
On vérifie immédiatement que sous les  hypothèses suivantes :
\ben
\item $f$ ne s'annule pas est $C^1$-limitée et $\inf(f(x)) \gnsim 0$
\item $\varphi > 0 $ et $\psi >0$
\item $z \gnsim 1 \Longrightarrow \psi(z) \;\mathrm{et}\;$ {i.g.}
\item $z \lnsim 1 \Longrightarrow \varphi(z)\sim 0  \;\mathrm{et}\; \psi(z)\sim 0$
\fen
on a :
\bitbul 
\item$z \lnsim 1 \Longrightarrow   z\sqrt{1+\varphi (z)  } ( f(x)-\psi(z)) \sim zf(x)   \quad \quad \quad z \gnsim 1 \Longrightarrow \frac{dz}{dt} = - \mathrm{i.g.}$  
\fit 
Donc, en dessous de la droite $z = 1$ on a $\Sigma \sim \tilde{\Sigma}$\\\\
\textbf{Le cas attractif :  $f>0$}. Le même argument que dans le cas lent-rapide (appendice \ref{halocourbelente}) permet d'affirmer que si $(x(t),y(t))$ est la trajectoire issue du point $a$ (au dessus de la droite $z =1$) il existe $t \sim0$ tel que $z(t) \sim 1$.\\
\begin{minipage}{0.5 \textwidth}
Partons maintenant du point $c$ situé en dessous de $z = 1$. Soit $(\tilde{x}(t),\tilde{y}(t)) $ la solution de $\tilde{\Sigma}$ à l'instant $0$ et $\tau_1$ l'instant où $\tilde{y}(\tau_1) = 0$.  Puisque $\inf(f(x)) \gnsim 0$ pour tout $\alpha \gnsim 0$ on a
$$\tau  \leq \alpha\,\longrightarrow \tilde{y}(\tau_1-\tau) \lnsim 1$$
et donc  $\Sigma \sim \tilde{\Sigma}$ ; d'après le lemme de Gronwall (appendice \ref{gronwall})  on a donc :
$$ 0 \lnsim \alpha \leq \tau_1  \Longrightarrow |y(t-\alpha)-\tilde{y}(t-\alpha)| \sim 0$$
L'ensemble :
\end{minipage}
\begin{minipage}{0.5 \textwidth}
\includegraphics[width=1\textwidth]{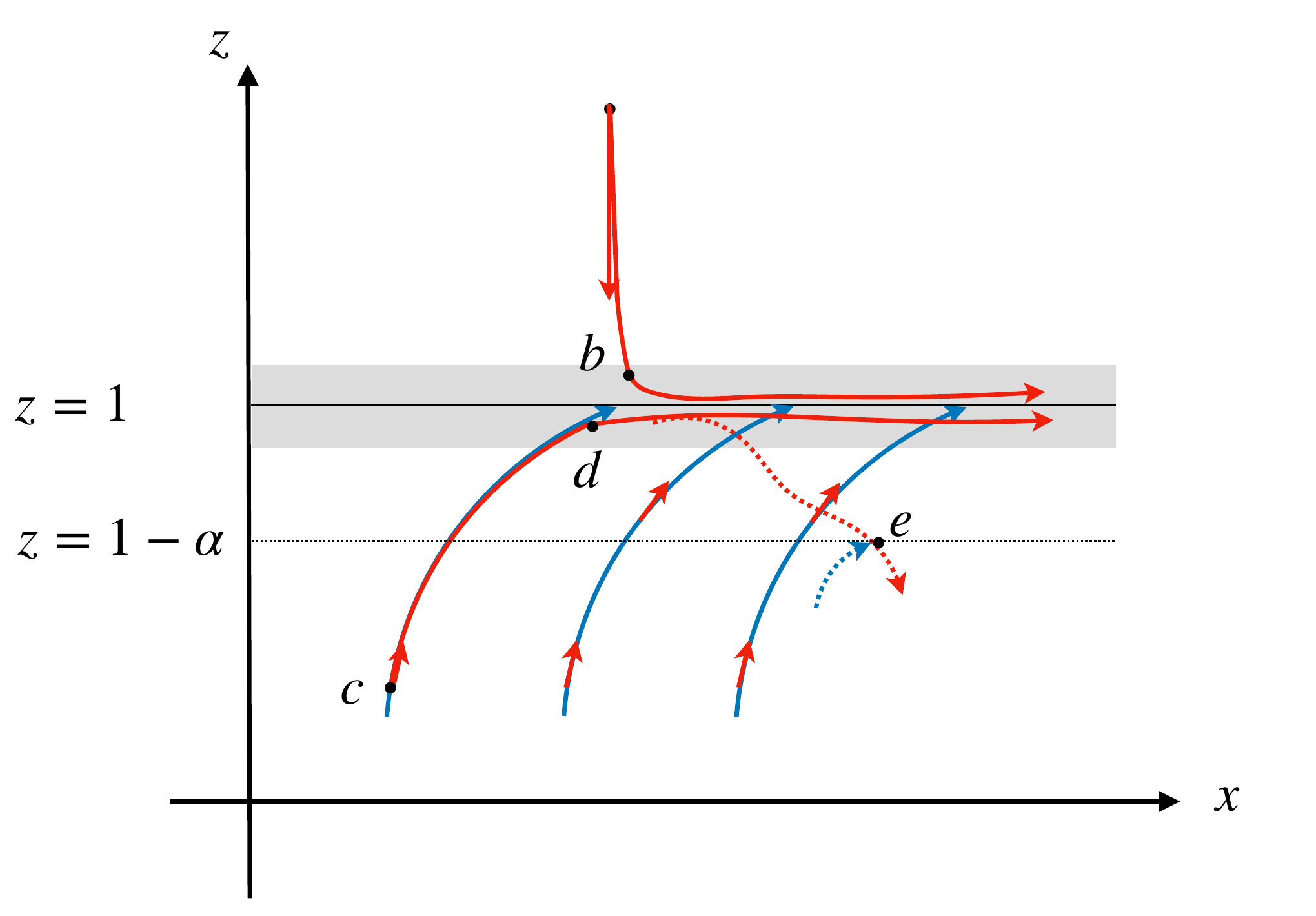}
\end{minipage} \\[6pt]
$$E = \left \{\alpha \in [0,\tau_1] : \frac{ |y(t-\alpha)-\tilde{y}(t-\alpha)|}{\alpha} < 1 \right \} $$
contient tous les $\alpha$ non infiniment petits, il contient donc un $\alpha \sim 0$ ce qui montre qu'il existe $t_1\sim \tau _1$ tel que $(x(t_1),y(t_1)) \sim (\tilde{x}(t_1),\tilde{y}(t_1)) $. La trajectoire de $\Sigma$  issue de $c$ à l'instant $0$ pénètre dans le halo de $z = 1$ en un instant $t_1$ infiniment proche de l'instant $\tau_1$ où $(\tilde{x}(t),\tilde{y}(t)) $ atteint $z = 1$ et reste infiniment proche de $(\tilde{x}(t),\tilde{y}(t)) $ sur $[0, t_1]$. A partir de cet instant la trajectoire reste infiniment proche de $z = 1$ ; en effet pour tout $\alpha \nsim 0$ si $t_2$ est le premier instant où $y(t) = 1-\alpha$ au point $e = (x(t_2),y(t_2)$ on a une contradiction avec le fait qu'en ce point le champ pointe strictement vers le haut.\\\\
\textbf{Le cas traversant :  $f< 0$}.\\
\begin{minipage}{0.5 \textwidth}
Soit le point $a =(x_0,y_0)$ au dessus de $z = 1$, le point $ c= (x_0,1)$ et $(\tilde{x}(t),\tilde{y}(t))$  la trajectoire de $\tilde{\Sigma}$. La trajectoire de $\Sigma$ issue de $a$ à l'instant $0$, pénètre  dans le halo de $z = 1$  en un instant $t_1 \sim 0$, en ressort pa dessous en un instant $t_2 \sim 0$ puis reste infiniment proche de la trajectoire $(\tilde{x}(t),\tilde{y}(t))$. La  durée de la traversée du halo de $z = 1$ est infiniment petite  parce que $\frac{dz}{dt} < z f(x) \leq \min(zf(x)) \lnsim 0$.
\end{minipage} $\quad$
\begin{minipage}{0.45 \textwidth}
\includegraphics[width=1\textwidth]{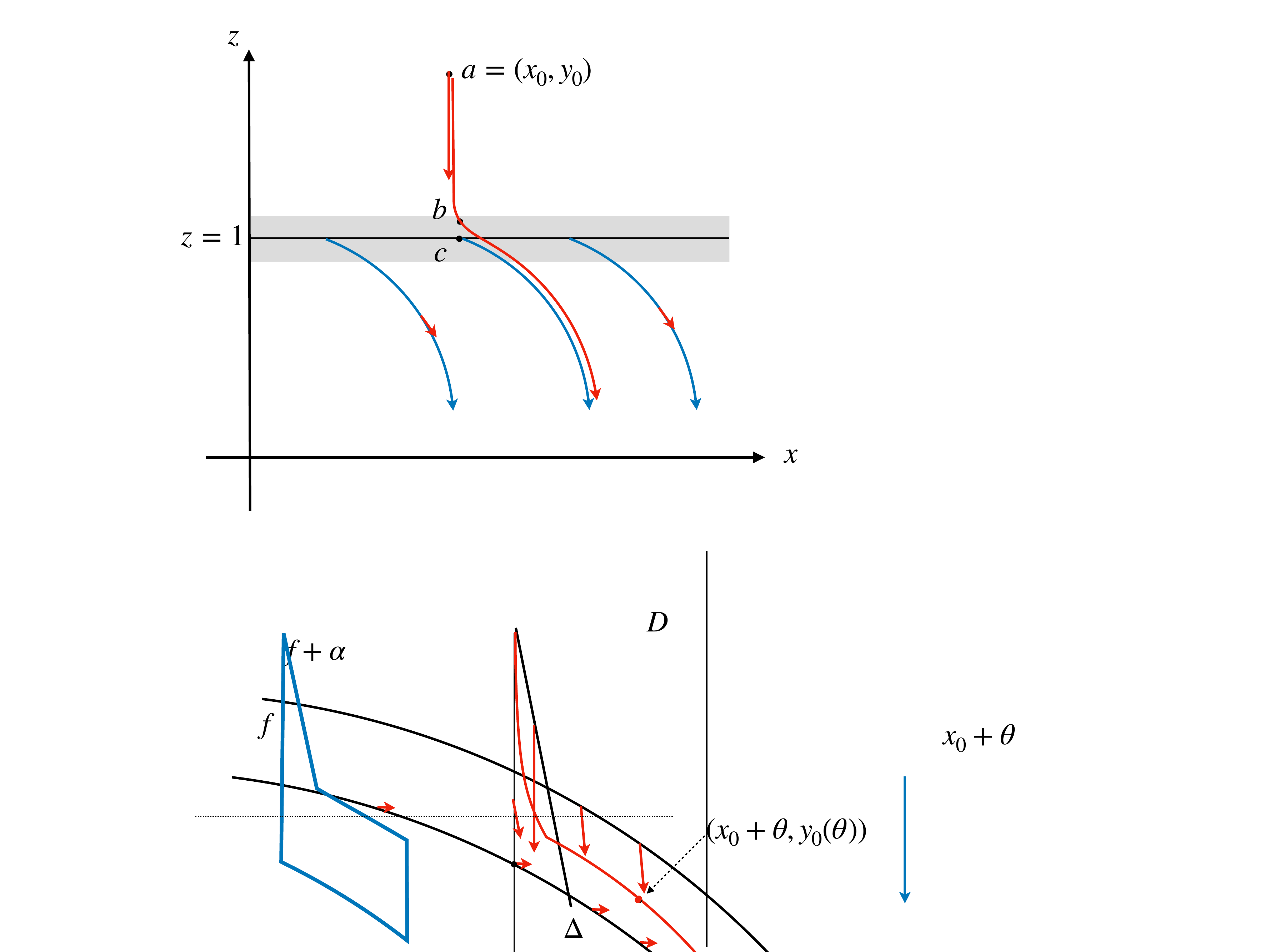}
\end{minipage} \\[6pt]

\subsection{Passage d'un changement de signe de $f$.}\label{changementdesigne}
  On considère pour $z > 0,\; \rho \lnsim 0, \;\eps \sim 0$ le système :
 \beq \label{eqchsigne}
F \quad \quad \quad  \left\{
\begin{array}{lcl}
\displaystyle \frac{dx}{dt}& =&1 \\[6pt]
\displaystyle \frac{dz}{dt} & =&\displaystyle z\sqrt{  1+\exp \left(2 \frac{\rho-\ln(z)}{\eps} \right)    } \left( f(x)-z^{\frac{1}{\eps}}\right)
\end{array} 
\right.
\feq
 et l'on fait les hypothèses suivantes :
 \begin{hypothese} Hypothèses sur $f$:
 \bitbul
\item Pour $ x > 0$ la fonction $f$ est égale à une fonction $f^+$, $C^1$-limitée  définie sur $\Rmat$ strictement négative pour $x>0$
\item Pour $ x < 0$ la fonction $f$ est égale à une fonction  $f^-$,  $C^1$-limitée  définie sur $\Rmat$ strictement positive pour $x<0$
 \fit
 \end{hypothese}
On introduit le  champ de vecteurs définie par  :
  \beq \label{eqchsignealphaplus}
\tilde{F}  \quad \quad \quad  \left\{
\begin{array}{lcl}
\displaystyle \frac{dx}{dt}& =&1 \\[6pt]
\displaystyle \frac{dz}{dt} & =&\displaystyle z f(z)
\end{array} 
\right.
\feq
 et l'on note respectivement   $\big(x(t,(x_0,z_0)),z(t,(x_0,z_0))\big)$ et $\big(\tilde{x}(t,(x_0,z_0)),\tilde{z} (t,(x_0,z_0))\big)$ les solutions de condition initiale $(x_0,z_0)$  à l'instant $0$ des champs  $F$ et $\tilde{F}$.
 \begin{lemme}
Soit $(x_0,z_0)$  une condition initiale telle que $x_0 \sim 0$ et $z_0 \sim 1$. Pour tout $t$ limité et tant que $\emat^{\rho} \leq \big(\tilde{x}(t,(0,1)),\tilde{z}(t,(0,1))\big)\leq 1$  :
$$ \big(x(t,(x_0,z_0)),z(t,(x_0,z_0))\big ) \sim  \big(x(t,(0,1)),z(t,(0,1))\big)$$
\end{lemme}
 \textbf{Démonstration} : On introduit les deux familles de champs de vecteurs définies pour $\alpha >0$ par :
 
  \beq \label{eqchsignealphaplusmoins}
F^{+\alpha}  \quad  \left\{
\begin{array}{lcl}
\displaystyle \frac{dx}{dt}& =&1 \\[6pt]
\displaystyle \frac{dz}{dt} & =&\displaystyle (1+\alpha)zf(z)
\end{array} 
\right.
\quad \quad
F^{-\alpha}  \quad   \left\{
\begin{array}{lcl}
\displaystyle \frac{dx}{dt}& =&1 \\[6pt]
\displaystyle \frac{dz}{dt} & =&\displaystyle (1-\alpha)zf(z)
\end{array} 
\right.
\feq
Sur la figure ci-contre les segments $\arc{b,c}$, $\arc{c,d}$  et $\arc{a,e}$ qui définissent le contour fermé $\partial \Gamma_{\alpha} =\arc{a,b,c,d,e,a}$ qui entoure le domaine $\Gamma_{\alpha}$   sont respectivement :\\
 \begin{minipage}{0.50 \textwidth} 
\bitbul
\item$\arc{b,c}$  : La portion de trajectoire de $F^{-\alpha}$ comprise entre $x = -\alpha$ et $x = 0$  qui passe par le point $(0,1-\alpha)$.
\item $\arc{c,d}$ : La portion de trajectoire de $F^{+\alpha}$ comprise entre $x = 0$ et $x = t $  qui passe par le point $(0,1-\alpha)$.
\item $\arc{a,e}$ La portion de trajectoire de $F^{-\alpha}$ comprise entre $x = -\alpha $ et $x = t $  qui passe par le point $(0,1+\alpha)$.
\fit
La trajectoire rouge est la solution de $F$ issue d'une condition initiale infiniment proche de $(0,1)$ et la trajectoire verte la solution de $\tilde{F}$ issue de $(0,1)$.
 \end{minipage}$\quad$
  \begin{minipage}{0.50 \textwidth}
 \includegraphics[width=1\textwidth]{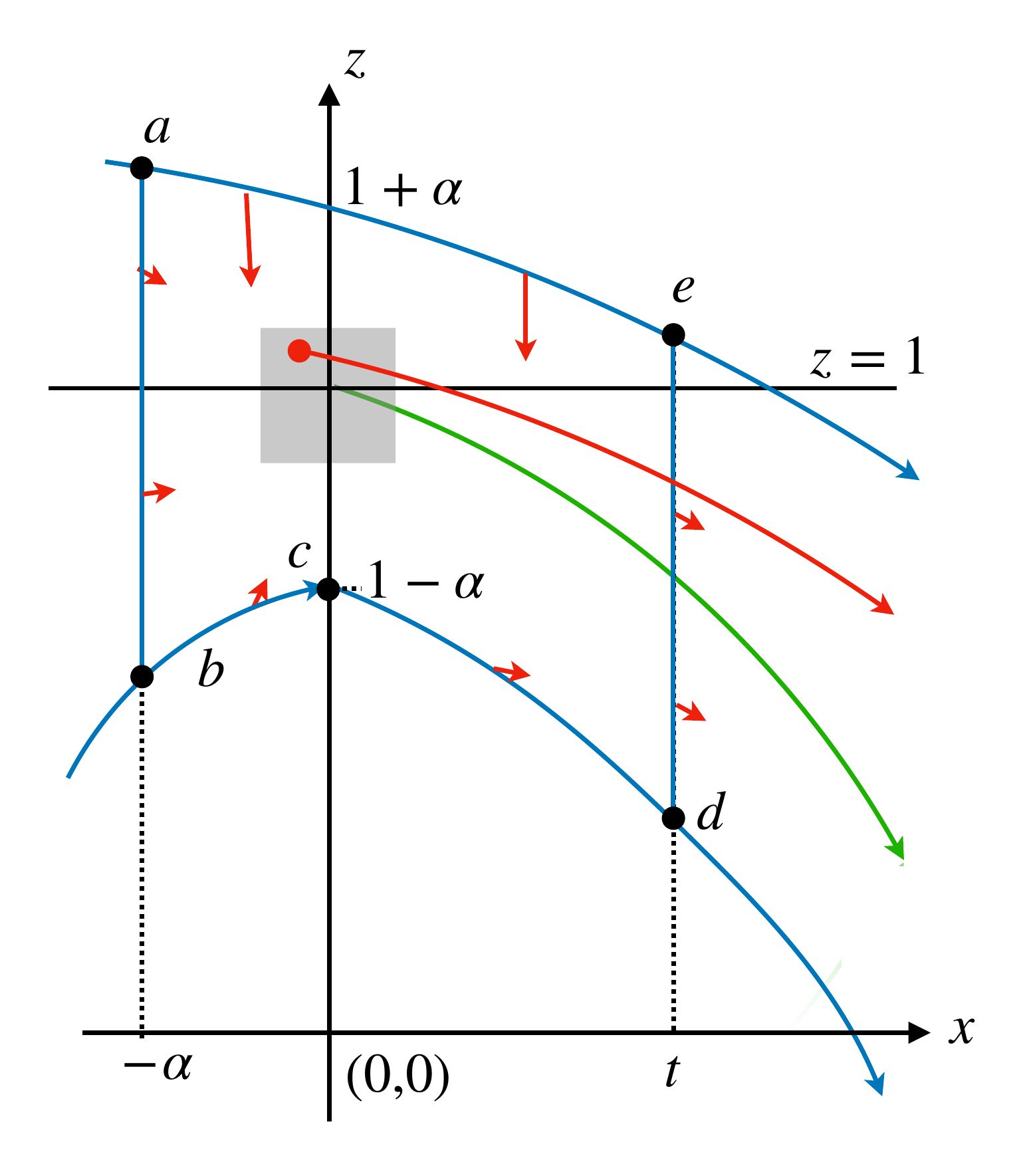}
 \end{minipage}
On vérifie aisément que si $\alpha \gnsim 0$ le champ  champ $F$ pointe strictement vers l'intérieur de $\Gamma_{\alpha}$ sur tous les segments du contour $\partial \Gamma_{\alpha} $ sauf sur $\arc{d,e}$ ; d'autre part, toujours si $\alpha \gnsim 0$ le point $(x_0,z_0) \sim (0,1)$ appartient à $\Gamma_{\alpha}$  et par suite (c.f. le théorème \ref{piege}) la trajectoire de $F$ qui en est issue 
sort de $\Gamma_{\alpha}$ en un point du segment $\arc{d,e}$. Comme les champ $F^{-\alpha}$ et $F^{+\alpha}$ tendent vers $F$ et, compte tenu des 
hypothèses de régularité sur $f$, quand $\alpha$ tend vers $0$ la longueur du segment $\arc{d,e}$ tend vers $0$ et par suite $ |z(t,(x_0,z_0))-\tilde{z}(t,(0,1))|\leq \delta$ pour tout $\delta \gnsim 0$ ce qui achève la démonstration.\\
$\Box$


\fin